\documentclass{amsart}
\usepackage{tikz}
\usetikzlibrary{arrows,matrix}
\usepackage[all]{xy}

\begin{document}
\title[Arithmetic differential geometry in the arithmetic PDE setting, I]
{Arithmetic differential geometry in the arithmetic  PDE setting, I: connections}
\author{Alexandru Buium}
\address{Department of Mathematics and Statistics,
University of New Mexico, Albuquerque, NM 87131, USA}
\email{buium@math.unm.edu} 

\author{Lance Edward Miller}
\address{Department of Mathematical Sciences,  309 SCEN,
University of Arkansas, 
Fayetteville, AR 72701}
\email{lem016@uark.edu}

\def \cH{\mathcal H}
\def \cB{\mathcal B}
\def \d{\delta}
\def \ra{\rightarrow}
\def \bZ{{\mathbb Z}}
\def \cO{{\mathcal O}}
\newcommand{\Hom}{\operatorname{Hom}}

\def\uy{\underline{y}}
\def\uT{\underline{T}}

\newcommand{\OCp}{ {\mathbb C}_p^{\circ}  }

\newcommand{\sbt}{\,\begin{picture}(-1,1)(-1,-2)\circle*{2}\end{picture}\ }

\newtheorem{THM}{{\!}}[section]
\newtheorem{THMX}{{\!}}
\renewcommand{\theTHMX}{}
\newtheorem{theorem}{Theorem}[section]
\newtheorem{corollary}[theorem]{Corollary}
\newtheorem{lemma}[theorem]{Lemma}
\newtheorem{proposition}[theorem]{Proposition}
\theoremstyle{definition}
\newtheorem{definition}[theorem]{Definition}
\theoremstyle{definition}
\newtheorem{convention}[theorem]{Convention}
\newtheorem{notation}[theorem]{Notation}
\newtheorem{remark}[theorem]{Remark}
\newtheorem{example}[theorem]{\bf Example}
\numberwithin{equation}{section}
\subjclass[2000]{11S31, 17B56, 53B05, 53C21}

\maketitle


\begin{abstract}
This is the first in a series on papers developing an
arithmetic PDE analogue  of Riemannian geometry.
The role of partial derivatives is played by Fermat quotient operations with respect to
several Frobenius elements in the absolute Galois group of a $p$-adic field. Existence and uniqueness of geodesics and of  Levi-Civita and Chern connections  are proved in this context.  In a sequel to this paper a theory of arithmetic  Riemannian curvature 
and  characteristic classes will be developed.
\end{abstract}

\section{Introduction}

The present paper is Part 1 of a series of papers;  its sequel \cite{BMadg2} will be referred to as Part 2. The aim of this series is to develop an arithmetic PDE version of the arithmetic ODE theory  in \cite{BD, Bu17, Bu19}. 

In \cite{BD, Bu17, Bu19} an arithmetic analogue of differential geometry was developed in which the role of functions is played by integers in $p$-adic fields and the role of 
derivations was played by $p$-derivations described by Fermat quotients where $p$ is a prime integer. Analogues of connections and curvature were introduced and analogues of familiar connections, for example Chern connections and Levi-Civita connection,  were studied.  These connections turned out to be {\it not} flat,  i.e., the integers have some `inherent curvature.' 
For a fixed prime $p$ this theory was developed in loc.cit in an unramified setting, over the completed valuation ring of the maximal unramified extension of $\mathbb Q_p$. A drawback, however, is that this ring has  a unique $p$-derivation, and hence the corresponding metrics/connections should be viewed as inherently of `cohomogeneity one' which made the theory  essentially an arithmetic ODE theory. 

The aim of the present paper and its sequel is to enhance the arithmetic ODE theory in \cite{BD, Bu17, Bu19} to an arithmetic PDE theory. Our approach involves a complete reformulation of the foundations and  admits a much wider natural scope. It not only resolves the cohomgeneity one limitation, but also  resolves a fundamental problem described in \cite{Bu17}; specifically we are able  to introduce here a theory of arithmetic geodesics not present in the ODE formulation.

The critical idea that allows for this framework to move past the ODE setting, is precisely that one must {\it not} restrict considerations to an unramified context. Indeed, the consideration of the (arbitrarily) ramified case allows one to implement the framework of purely arithmetic PDEs as introduced in \cite{BMP}. 

 The present paper is essentially self contained and can be viewed as a  natural next step for the theory initiated in \cite{Bu95} and developed, in particular, in \cite{Bu05,Bu17,BMP},
as explained in what follows. On the other hand, all of this work can be viewed as an instance of the general effort to discover new analogies between numbers and functions; see, in particular, Manin's paper \cite{Man13}.

\subsection{Acknowledgements} The present work, as well as much of the previous work by the first author, were deeply influenced by Yu.I.Manin's  ideas on the analogies between numbers and functions. The first author was partially supported by the Simons Foundation through award 615356.

\subsection{Overview of  previous work} 
The first step in this circle of ideas was taken in \cite{Bu95} where
a theory of arithmetic ODEs was initiated. In \cite{Bu95} derivation operators were replaced by $p$-derivation operators $\delta_p$ acting on rings of integers of local fields. This theory led to  a series of Diophantine applications; see \cite{Bu95, Bu96, Bu97, Bu05, BP09}. 
In particular, in   \cite{Bu95, Bu97}, analogues of the classical Manin maps \cite{Man63}  were constructed and   arithmetic analogues of Manin's Theorem of the Kernel were  proved. On the other hand this theory provided some striking arithmetic analogues of Hamiltonian systems appearing in mathematical physics \cite{BMan, BP17, BP18, Bu20}.

In \cite{Bu17}, this theory was adapted to give an arithmetic version of differential geometry. Among other things certain ODE versions of the PDEs appearing in classical Riemannian  geometry related to Chern and Levi-Civita  connections were developed. In particular, the arithmetic Christoffel symbols were shown to enjoy a surprising relationship with, and could  be viewed as matrix analogues of, Legendre symbols. 

The framework of \cite{Bu05} for arithmetic ODEs considers only `unramified solutions' i.e., solutions with coordinates in the completion $R$ of the maximum unramified extension of $\mathbb Z_p$. Subsequently, the $\delta$-overconvergence machinery in \cite{BS11, BM20} allowed one to add `ramified solutions' to the main arithmetic ODEs of the theory, i.e., solutions  in the ring of integers $R^{\text{alg}}$ of the algebraic closure $K^{\text{alg}}$ of $K:=R[1/p]$. Sometimes even solutions in the ring of integers  of the complex $p$-adic field $\mathbb C_p$ can be attained. Similarly, the arithmetic differential geometry described in \cite{Bu17,Bu19} is restricted to unramified solutions and the present article and its sequel introduce in this theory ramified solutions and purely arithmetic PDEs. 

To contrast, a  theory of  arithmetic PDEs with $2$ `directions' was developed in \cite{BS10a, BS10b}. One of these directions was arithmetic and the other geometric. The solutions were  again  unramified, lying in the ring of power series $R[[q]]$ in one variable $q$. In  \cite{BMP} we introduced a   `purely arithmetic' variant of this idea namely, a PDE theory in which all the directions are `arithmetic,' and in which one considers solutions in $R^{\text{alg}}$. Here, surprising novel versions of Manin maps where constructed leading to enhancements of the Diophantine applications, new reciprocity laws, and more. 


We explain now a bit more what we mean by `purely arithmetic' PDEs. The starting point in \cite{BMP} is that one should envision not one but several arithmetic directions at every given prime. This comes from the fact that the absolute Galois group $\mathbb Q_p$ has several topological generators. One then builds the theory starting from an arbitrary finite collection  $\phi_1,\ldots,\phi_n$ of   Frobenius automorphisms of $K^{\textup{alg}}$. Remarkably this approach, combined with the $\delta$-overconvergence technique in \cite{BS11, BM20}, permits consideration of solutions to equations in $R^{\text{alg}}$. 

The main idea in \cite[Sec. 4.3]{Bu17} was to associate to every $n\times n$ symmetric matrix $q$ with coefficients in $R$, viewed as an analogue of a metric, 
an $n$-tuple of $p$-derivations on 
the ring of functions on $\textup{GL}_n$ over $R$. The $p$-derivations were subject to  an appropriate symmetry condition and could be viewed as an analogue of the Levi-Civita connection with respect to `one arithmetic direction.' This further led to an analogue of Riemannian curvature. In \cite{Bu19} a global version of this was developed: 
 for every rational prime $p$ we considered $n$ `arithmetic directions' given by $n$ primes  dividing $p$ in a given number field and we considered {\it one} `metric' with coefficients in the number field. However, this setting  was `reducible' to the case of {\it one} `arithmetic direction'  and $n$ `metrics,' making the theory, again, an ODE theory. Also, in this setting, only unramified solutions were allowed. 
 
 In the present paper, we utilize the framework in \cite{BMP} to put forward a `genuinely PDE' formalism that corresponds to the case of `$n$ arithmetic directions' at one prime and one `metric' $q$ specifically to construct arithmetic Levi-Civita and Chern connections. Part 2 of this paper will  address questions about  curvature and characteristic classes.

We summarize our comparison between various approaches to arithmetic Riemannian geometry  in the following table; $N_{\textup{pr}}$ is the number of primes, $N_{\textup{ari}}$ is the number of arithmetic directions and $N_{\textup{metric}}$ is the number of  metrics: 

\bigskip

\begin{center}
\begin{tabular}{||c|c|c|c|c||} \hline \hline
reference  & $N_{\textup{pr}}$ & $N_{\textup{ari}}$ & $N_{\textup{metric}}$  & ramified solutions defined\\ 
\hline
\cite{Bu17} & 1 & 1 & 1  & NO \\
\hline
\cite{Bu19} & 1 & 1 & n  & NO\\
\hline
\textup{this paper} & 1 & n & 1  & YES\\
\hline
\end{tabular}
\end{center}

\bigskip

We end this review by contrasting two analogies present in $p$-adic theories; let us  call them
 here the {\it number-function} analogy and the {\it number-number} analogy. 
 In the number-function analogy $p$-adic {\it numbers} are viewed as analogous to 
{\it functions of one variable}. In the number-number analogy $p$-adic {\it numbers} are viewed as analogous to {\it real numbers} (rather than functions). 
The classical `$p$-adic arithmetic' (from Hensel's Lemma to the ramification and Galois theory of $p$-adic fields)
can be viewed as an instance of the number-function analogy. In addition, 
   the theory of arithmetic differential equations in \cite{Bu05} and the arithmetic differential geometry in \cite{Bu17} are based on this same number-function analogy.
In contrast with this the field of  `$p$-adic analysis' (including $p$-adic measure theory \cite[Ch. 2]{Ko84}, $p$-adic analytic function theory \cite[Ch. 4]{Ko84},
the theory of $p$-adic differential equations  \cite{Ked}, `$p$-adic physics'  \cite{VVZ}, etc.) 
is based on the number-number analogy. There are interesting situations, however, when 
the mathematical objects emerging from  these two different analogies can be made to
interact. For instance,  $p$-adic $L$-functions (which are built via a number-number analogy) play a role in Iwasawa theory (which fits into the number-function analogy); see \cite{Wa}. 
As another example certain $p$-adic analytic functions considered by Dwork  \cite[Ch.V]{Ko84} (which fit into the number-number analogy) were used by him to prove the rationality of zeta functions of varieties over finite fields (which fits into the number-function analogy). In a similar vein
it is conceivable that interactions exist between \cite{Bu05,Bu17} on the one hand and 
\cite{Ked,VVZ} on the other;  we do not pursue this here.

\subsection{Framework  and results of this paper} 
We sketch here some of our main concepts and results.  Let  $\pi$ be a prime element in a finite Galois extension of $\mathbb Q_p$. We denote by $R$  the completion of the maximum unramified extension of $\mathbb Z_p$. Let $R_{\pi}:=R[\pi]$, $k:=R_{\pi}/\pi R_{\pi}=R/pR$, $K_{\pi}:=K(\pi)$. 
The degree $[K_{\pi}:K]$ of the field extension $K_{\pi}/K$ equals the ramification index $e$ of this extension. The ring $R_{\pi}$ will be viewed as an arithmetic analogue of the ring of smooth functions on a smooth manifold $M$. Fix an integer $s\in \mathbb N$.
 In what follows by  a {\bf higher $\pi$-Frobenius lift} of degree $s$ on a  ring $S$ containing $R_{\pi}$ we understand  a ring endomorphism $\phi$ of $S$ reducing to the $p^s$-power Frobenius modulo $\pi$. We start with an $n$-tuple of $\pi$-Frobenius lifts, $(\phi_1^{(s)},\ldots,\phi_n^{(s)})$, of degree $s$ on $R_{\pi}$ and we consider the  {\bf higher $\pi$-derivations} of degree $s$, $\delta_i^{(s)}:R_{\pi}\rightarrow R_{\pi}$, 
 attached to $\phi_i^{(s)}$, defined by the formula 
 $$\delta_i^{(s)}(\lambda):=
 \frac{\phi_i^{(s)}(\lambda)-\lambda^{p^s}}{\pi},\ \ \ \lambda\in R_{\pi}.$$
 We view $\delta_1^{(s)},\ldots,\delta^{(s)}_n$  as an arithmetic analogue of a basis $X_1,\ldots,X_n$ for the module of vector fields on a parallelizable smooth manifold $M$. 
 We then fix an integer $N\geq 1$ and we consider the group scheme $G$ and the ring $\mathcal A$ defined by
 $$G:=\textup{GL}_{n/R_{\pi}}:=\textup{Spec}
 \ R_{\pi}[x,\det(x)^{-1}],\ \ \ \mathcal A:=\widehat{\mathcal O(G)}=R_{\pi}[x,\det(x)^{-1}]^{\widehat{\ }}$$  
   where we denote by $x=(x_{ij})$  an $N\times N$ matrix of indeterminates, the $\widehat{\ }$ sign means $p$-adic completion, and $\mathcal O(G)$ is the global ring of functions for $G$.  The ring $\mathcal A$ is viewed as an arithmetic analogue of the ring of smooth functions on the principal  bundle attached to a vector bundle of rank $N$ over the manifold $M$.
 We define a {\bf $\pi$-connection} on $G$ of degree $s$ to be an $n$-tuple of higher $\pi$-Frobenius lifts $((\phi_1^{(s)})^G,\ldots,(\phi_n^{(s)})^G)$ of degree $s$ on 
 $\mathcal A$ extending those on $R_{\pi}$. To give a $\pi$-connection is equivalent to giving the $n$-tuple of {\bf higher $\pi$-derivations} 
$(\delta_i^{(s)})^G:\mathcal A\rightarrow \mathcal A$  of degree $s$ attached to $(\phi_i^{(s)})^G$, defined by the formula 
 $$(\delta_i^{(s)})^G(a):=\frac{(\phi_i^{(s)})^G(a)-a^{p^s}}{\pi},\ \ \ a\in \mathcal A.$$
 The tuple  $(\delta_1^{(s)})^G,\ldots,(\delta_n^{(s)})^G$ is analogous to a usual connection on a principal  bundle attached to a vector bundle. The numbers
 \begin{equation}\label{4numbers}
 e,s,n,N\end{equation}
 introduced above are not assumed to be related in general.
 In case $e=n=N$ one can introduce {\bf parallel transport} and {\bf geodesics}. We will not review these concepts in this Introduction; suffices to say that  existence and uniqueness results will be proved for parallel transport and geodesics;
 see Theorems \ref{eupar} and \ref{eug}. 
 
 Reverting to the case where 
 the numbers (\ref{4numbers}) are unrelated
 the next concept to be considered will be that of {\bf metric} by which we will understand a
 symmetric matrix $q\in \textup{GL}_N(R_{\pi})$. This is the analogue of a Riemannian  metric on a vector bundle on a  manifold.
 In Part 2 we will consider metrics for each degree $s$ so it is useful to we write, in what follows, $q=q^{(s)}$. 
  We will say that  a $\pi$-connection as above is {\bf metric} with respect to $q^{(s)}$
   if the following diagrams are commutative for all $i$:
 $$
 \begin{array}{rcl}
\mathcal A& \stackrel{(\phi_{i,0}^{(s)})^G}{\longrightarrow} & \mathcal A \\
 \cH_{q^{(s)}} \downarrow &\ &\downarrow \cH_{q^{(s)}}\\
 \mathcal A & \stackrel{(\phi_i^{(s)})^G}{\longrightarrow} & \mathcal A\end{array}
$$
where $(\phi_{i,0}^{(s)})^G$ are the $\pi$-Frobenius lifts that send the matrix $x=(x_{jk})$ into the matrix $x^{(p^s)}:=(x_{jk}^{p^s})$ and $\mathcal H_{q^{(s)}}$ is the $R_{\pi}$-algebra homomorphism defined by $\mathcal H_{q^{(s)}}(x)=x^t q^{(s)} x$. Although this is not immediately visible, the above definition is an arithmetic analogue of the usual metric condition for connections on the principal  bundle attached to a vector bundle; cf. \cite[pg. 64--65]{Bu17} or the Appendix of Part 2 of this paper, \cite{BMadg2}. 

We next assume the numbers $n,N$ in (\ref{4numbers}) are subject to the condition $N=n$.
We define a {\bf torsion symbol} to be 
  an $n$-tuple $L^{(s)}=(L^{k(s)})$ of $n\times n$ 
  antisymmetric matrices $L^{k(s)}=(L^{k(s)}_{ij})$, 
  with entries in the ring $\widehat{R_{\pi}[y]}$ where $y$ is an $n$-tuple of $n\times n$ matrix indeterminates. 
    In the simplest case of our theory, to be fully explored in Part 2 of this paper,  the matrices $L^{k(s)}$ above will actually have entries in the set $\{-1,0,1\}$ and will reflect the commutation relations between the $\phi_i^{(s)}$s. If $M$ is a Lie group and $X_1,\ldots,X_n$ is a basis for the space of left  invariant vector fields
  then the elements $L^{k(s)}_{ij}$ are arithmetic analogues of the structure constants of the Lie algebra of $M$.
  Define the {\bf Christoffel symbols} of the second kind $\Gamma_{ij}^{k(s)}\in \mathcal A$ as the $jk$-entries of the matrices $\Gamma_i^{(s)}$ whose transposed are given by:
  $$
\Gamma_i^{(s)t}:=(x^{(p^s)})^{-1} (\delta_i^{(s)})^G x \in \textup{Mat}_n(\mathcal A).$$ 
This definition is analogous to the classical definition of Christoffel symbols; cf. \cite[p. 44]{Bu17}.
 We also consider the matrices
    \begin{equation}
\Lambda_i^{(s)}:=(x^{(p^s)})^{-1}(\phi_i^{(s)})^G (x)\in \textup{GL}_N(\mathcal A)\end{equation}
and we  set $\Lambda^{(s)}=(\Lambda_1^{(s)},\ldots,\Lambda_n^{(s)})$.

From these choices, we will combine the torsion symbols with $\Lambda^{(s)}$ to obtain elements $L^{k(s)}_{ij}(\Lambda^{(s)}) \in \mathcal A$ which measure symmetry among the Christoffel symbols in the following sense. We say  that a $\pi$-connection 
   of degree $s$ is {\bf symmetric} with respect to $L^{(s)}$
   if 
    \begin{equation}
\label{uf2}
\Gamma^{k(s)}_{ij}-\Gamma^{k(s)}_{ji}= L^{k(s)}_{ij}(\Lambda^{(s)}).\end{equation}
This is an arithmetic analogue of the symmetry condition for  connections on the frame bundle of the  manifold $M$; cf. the Appendix to Part 2.  One of our main results, Theorem \ref{LCC},
 will imply the following arithmetic analogue of the Fundamental Theorem of Riemannian Geometry.
  
 \begin{theorem}\label{boooster}
 For every metric $q^{(s)}$ and every torsion symbol $L^{(s)}$ there exists a unique $\pi$-connection of degree $s$ that is  metric with respect to $q^{(s)}$ and symmetric with respect to $L^{(s)}$.
 \end{theorem}

This $\pi$-connection will be referred to as the  {\bf arithmetic  Levi-Civita connection} of degree $s$ attached to $q^{(s)}$. As in the ODE case its Christoffel symbols are related to the Legendre symbol
and can be  viewed as a matrix analogue of the latter. For a discussion of the connection with  the Legendre symbol we refer to Part 2. 

We will develop a similar theory  for the {\bf Chern connection} in which $N$ is not necessarily equal to $n$; for the existence and uniqueness of Chern connection  see Theorem \ref{ch}. 

We will close our paper by considering what happens if we vary $\pi$; the main result here is Theorem \ref{javasp} which states that our Levi-Civita and Chern connections enjoy a {\bf $\delta$-overconvergence} property in the sense of \cite{BMP}; this morally allows one to formulate the theory `over $K^{\textup{alg}}$.' 


\subsection{Notation and conventions}\label{generall}
In this paper and its Part 2 we adopt the following conventions.
Unless otherwise stated,  monoids will not necessarily have an identity. We denote by $\mathbb N$ the additive monoid of positive integers.
Associative $\mathbb Z$-algebras will   not be assumed commutative or with identity. By a ring we will mean a commutative associative $\mathbb Z$-algebra with identity. 
  For a (not necessarily algebraic) field extension $L/F$, denote by $\mathfrak G(L/F)$ the group of $F$-automorphisms of $L$.   The superscript `alg' for a field and `ur' for a local field will mean `algebraic closure' and `maximum unramified extension,' respectively. The superscript `t' for a matrix will mean `transpose.' 
We fix throughout an odd prime $p\in \mathbb Z$. For any ring $S$ and any Noetherian scheme $X$ we always denote by   $\widehat{S}$ and $\widehat{X}$ the respective $p$-adic completions.  For $X$  a scheme or a formal scheme 
we denote by $\mathcal O(X)$ 
 its global ring of functions. 

It is standard notation in \cite{Bu17}, for a matrix $a$  with entries $a_{ij}$ 
in a ring to denote by $a^{(p^d)}$ the matrix whose entries are $a_{ij}^{p^d}$. We adopt this with one notable exception. Our use of higher lifts also defines for a matrix $a$ and natural number $s \geq 1$ a matrix $a^{(s)}$ in a different sense for example in $\Gamma_i^{(s)}$, $\Lambda_i^{(s)}$, $q^{(s)}$, etc.,  where $(s)$ refers  to the degree. Thus when $s=p^d$ for some $d\geq 1$ there is a conflict between our use of the superscript $(s)$. To resolve this conflict we make the convention that  $a^{(p^d)}$ will mean $(a_{ij}^{p^d})$ if the exponent is explicitly written as a $p$-power, even if this power is $p^s$; which will occur. The exponent $(s)$ will be used to denote the degree if $s$ is {\it not} explicitly written as a $p$-power.

\subsection{Leitfaden} In Section 2 we review  the main arithmetic PDE setting in \cite{BMP} and we also provide some complements on the arithmetic ODEs of \cite{Bu95,Bu05}.  Section 3 develops the general  formalism of connections in the arithmetic PDE setting, with applications to geodesics and to  Levi-Civita and Chern connections. In Section 4 we offer an invariant look at our main concepts by placing them in the setting of arithmetic jet spaces and arithmetic Maurer-Cartan equations (logarithmic derivative map). 
In Section 5 we use the framework of Section 4 to address the compatibility of our constructions with varying $\pi$ and prove the $\delta$-overconvergence of our Chern and Levi-Civita connections. As already mentioned, Part 2 will deal with curvature and characteristic classes. We will also consider there Legendre symbols, gauge actions, torsors, and the construction of canonical metrics and torsion symbols from Galois theoretic data.

\section{Preliminaries}

We start with a simplified review of the main setting  in \cite{BMP,BM20} as well as introduce in detail notions of higher Frobenius lifts.

\subsection{Frobenius automorphisms}

Consider the field $\mathbb Q_p$ of $p$-adic numbers with absolute value $|\ |$ normalized by $|p|=p^{-1}$.
Fix an algebraic closure of ${\mathbb Q}_p$ denoted by ${\mathbb Q}_p^{\text{alg}}$, let
${\mathbb Q}_p^{\text{ur}}$  be the maximum unramified extension of ${\mathbb Q}_p$
inside ${\mathbb Q}_p^{\text{alg}}$, let $K$ be the metric completion of ${\mathbb Q}_p^{\text{ur}}$ and let $K^{\text{alg}}$ be the algebraic closure of $K$ in the metric completion ${\mathbb C}_p$
of ${\mathbb Q}_p^{\text{alg}}$. Recall that $K$ is not an algebraic extension of $\mathbb Q_p$.
We still denote by $|\ |$ the induced absolute value on all of these fields.  We denote by
${\mathbb Z}_p^{\text{ur}}, {\mathbb Z}_p^{\text{alg}}, R, R^{\text{alg}}$
the valuation rings of 
${\mathbb Q}_p^{\text{ur}}, {\mathbb Q}_p^{\text{alg}}, K, K^{\text{alg}}$, 
respectively.
We set
$R=\widehat{{\mathbb Z}_p^{\text{ur}}}$ and  $k:=R/pR\simeq \mathbb F_p^{\textup{alg}}$ the residue field\footnote{Throughout, there is a conflicting usage of $k$ as both the residue field and at times as an natural number index, e.g., $\Gamma_{jk}$ for the $(j,k)$-th entry in a matrix $\Gamma$. Context should make clear the intended use and no confusion should arise. }. 
  We denote by $\textup{Fr}:k\rightarrow k$ the $p$-power Frobenius and by $\textup{Fr}^s$ its $s$-th iterate, $s\in \mathbb N$.
 Recall that the natural ring homomorphism
  $\mathbb Q_p^{\text{alg}}\otimes_{\mathbb Q_p^{\text{ur}}} K\rightarrow K^{\text{alg}}$
  is an isomorphism, cf. \cite[Rmk. 2.2]{BMP}; surjectivity is a consequence of  Krasner's lemma.
  Recall also that the roots of unity in $R$ are exactly the roots of unity of  order prime to $p$ in $K^{\textup{alg}}$ and the canonical surjection $R\rightarrow k$ induces an isomorphism between the group  of roots of unity in $R$ and the group $k^{\times}$.

  \begin{definition}\label{lpdis}
  Let $L$ be a subfield of $\mathbb C_p$ containing $\mathbb Q_p$.
    A {\bf higher Frobenius automorphism} of degree $s\geq 1$ on  $L$  is a continuous automorphism  $\phi\in \mathfrak G(L/\mathbb Q_p)$  that  induces the $p^s$-power automorphism on the residue field of the valuation ring  of $L$. For any subfield $L_0$ of $L$ we denote by $\mathfrak F^{(s)}(L/L_0)$ the set of higher Frobenius automorphisms on $L$ of degree $s$ that are the identity on $L_0$. \end{definition}

  \begin{remark} \ 
  
  1) We usually write $\phi^{(s)}$ instead of $\phi$ for $\phi\in \mathfrak F^{(s)}(L/L_0)$ if we want to make the degree $s$ explicit. 
  
  2) For every $s\geq 1$ the set $\mathfrak F^{(s)}(K/\mathbb Q_p)$ has cardinality $1$ and the set $\mathfrak F^{(s)}(K^{\textup{alg}}/\mathbb Q_p)$ is a non-empty principal homogeneous space for the Galois group $\mathfrak G(K^{\textup{alg}}/K)$.
  
  3) A special case of what we call here `higher Frobenius automorphisms' are sometimes simply called `Frobenius automorphisms' in the literature, cf. e.g.,  \cite[p. 41]{N80}. Since our setting is more general than, for instance, that of \cite{N80}, and in order to align our exposition with that of \cite{BMP}, we prefer to use the slightly different terminology in Definition
  \ref{lpdis}.\end{remark}

 \begin{remark}
\label{oiap}
Before we continue our discussion recall  the following standard facts which follow from 
 \cite[p. 19, Props. 17 and 18]{Se79}. Let $F$ be a finite extension of $\mathbb Q_p$.
 
 1) If $\pi$ is a root of an Eisenstein polynomial with coefficients in the valuation ring $\mathcal O_F$ of $F$ and $E:=F(\pi)$  then $\mathcal O_E=\mathcal O_F[\pi]$ is the valuation ring of $F(\pi)$ and $\pi$ is a prime element of $\mathcal O_E$. 
 
 2) If  $E$ is a finite totally ramified extension of $F$ and $\pi$ is a prime element in $E$ (i.e., in  $\mathcal O_E$) then $\pi$ is a root  of an Eisenstein polynomial with coefficients in $\mathcal O_F$ and we have $\mathcal O_E=\mathcal O_F[\pi]$, hence $E=F(\pi)$.  If we assume in addition that $F\subset K$ and 
 we set $K_{\pi}:=K(\pi)$, 
$R_{\pi}:=R[\pi]$ then $R_{\pi}$ is the valuation ring of $K_{\pi}$ and $K_{\pi}$ is isomorphic to the tensor product $E\otimes_F K$.
\end{remark} 

\subsection{$\pi$-derivations}\label{subsectionbeloww}

As in \cite{BMP}, throughout this paper, we denote by $\Pi$ the set of all elements $\pi\in {\mathbb Q}_p^{\text{alg}}$ such that there is a finite Galois extension $E/\mathbb Q_p$
with the property that $\pi$ is a prime element of $E$. Recall from Remark \ref{oiap} that
if $F$ is the maximum unramified extension of $\mathbb Q_p$ contained in $E$ and $\pi$ is a prime element of $E$
then  $\pi$ is a root  of an Eisenstein polynomial with coefficients in $\mathcal O_F$ and we have $E=F(\pi)$.  
It follows  that ${\mathbb Q}_p^{\text{alg}}={\mathbb Q}_p^{\text{ur}}(\Pi)$. Also if $K_{\pi}=K(\pi)$ then 
$R_{\pi}:=R[\pi]$ is the valuation ring of $K_{\pi}$. Note that $K^{\text{alg}}=K(\Pi)$.
For $\pi,\pi'\in \Pi$ write $\pi'|\pi$ if and only if $K_{\pi}\subset K_{\pi'}$.
Note also that the set $\Pi$ consists exactly of the elements $\pi\in  {\mathbb Q}_p^{\text{alg}}$
which are roots of Eisenstein polynomials with coefficients in $\mathbb Z_p^{\textup{ur}}$
and for which the extension $\mathbb Q_p^{\textup{ur}}(\pi)/\mathbb Q_p$ is Galois.

\begin{remark}
We take the opportunity to correct here a typo in \cite{BM20}:  
 in the definition of $\Pi$ of Section 2.1 in loc.cit. the exponent `ur' in the condition 
  `$\mathbb Q_p^{\textup{ur}}(\pi)/\mathbb Q_p$ is Galois' was inadvertently omitted.  
\end{remark}

\begin{notation}
We will often use the following notation/convention. For every $\alpha\in k$ and $a\in R_{\pi}$ we write 
$$\alpha=a\ \ \textup{mod}\ \ \pi$$ if the image of $a$ in $k$ is $\alpha$.
\end{notation}

It will be convenient to utilize the following notion.

\begin{definition}
Let $A$ be an $R_{\pi}$-algebra.
By a {\bf  higher $\pi$-Frobenius lift} of degree $s$  for an $A$-algebra $\varphi\colon A\rightarrow B$
 we understand a ring homomorphism $\phi:A\rightarrow B$ such that the induced homomorphism $\overline{\phi}\colon A/\pi A\rightarrow B/\pi B$
equals the composition of the induced homomorphism $\overline{\varphi}:A/\pi A\rightarrow B/\pi B$ with the $p^s$-power Frobenius on $A/\pi A$. 
If $B=A$, and $\varphi=1_A$ we say that $\phi$ is a higher  $\pi$-Frobenius lift  on $A$ of degree $s$; if in addition $\pi=p$ and $s=1$ we simply say $\phi$ is a Frobenius lift.
In order to include the degree 
in the notation we usually write $\phi^{(s)}$ instead of $\phi$; if we further want to include $\pi$ in the notation we write $\phi^{(s)}_{\pi}$ instead of $\phi$.\end{definition}

Next, for every $s\in \mathbb N$, consider the polynomial
$$C_p^{(s)}(X,Y):=\frac{X^{p^s}+Y^{p^s}-(X+Y)^{p^s}}{p}\in \mathbb Z[X,Y].$$
Following  \cite{J85,Bu95} we introduce the following:

 \begin{definition}
 Let $A$ be an $R_{\pi}$-algebra and $B$ and $A$-algebra.
A map $\delta^{(s)}_{\pi}:A\rightarrow B$ is called a {\bf higher $\pi$-derivation} of degree $s$ if 
 $\delta^{(s)}(1)=0$ and
\[\begin{array}{rcl}
\d_{\pi}^{(s)}(x+y) & = &  \d^{(s)}_{\pi} x + \d^{(s)}_{\pi} y
+\frac{p}{\pi}\cdot C^{(s)}_p(x,y)\\
\ & \ & \ \\
\d_{\pi}^{(s)}(xy) & = & x^{p^s} \cdot \d^{(s)}_{\pi} y +y^{p^s} \cdot \d^{(s)}_{\pi} x
+\pi \cdot \d^{(s)}_{\pi} x \cdot \d^{(s)}_{\pi} y,
\end{array}\] for all $x,y \in A$.
The map $\phi_{\pi}^{(s)}:A\rightarrow B$, $\phi^{(s)}_{\pi}(x):=\varphi(x)^{p^s}+\pi\delta^{(s)}_{\pi} x$ is then a higher $\pi$-Frobenius lift of degree $s$. If $B$ is torsion free we get 
a bijection $\delta_{\pi}^{(s)}\mapsto \phi^{(s)}_{\pi}$ between the set of higher $\pi$-derivations of degree $s$ from $A$ to $B$ and the set of higher $\pi$-Frobenius lifts of degree $s$ from $A$ to $B$. 
We say that $\delta_{\pi}^{(s)}$ and $\phi_{\pi}^{(s)}$ are associated to each other. We sometimes identify elements $x\in A$ with the elements $\varphi(x)=x\cdot 1_B$. \end{definition}

\begin{remark}
Note the formula
$$\delta_{\pi}^{(s)} \pi\equiv \frac{\phi_{\pi}^{(s)}\pi}{\pi}\ \ \ \textup{mod}\ \ \pi\ \ \textup{in}\ \ R_{\pi}.$$
In particular, we have
$$\delta_{\pi}^{(s)} \pi\not\equiv 0\ \ \ \textup{mod}\ \ \pi.$$
Also note that if $A=B$ and $\phi^{(s)}_{\pi}\pi=\pi$ then $\phi^{(s)}_{\pi}$ and $\delta^{(s)}_{\pi}$ commute.
\end{remark}

\begin{definition}
A {\bf partial $\delta_{\pi}$-ring} of degree $s$ is an $R_{\pi}$-algebra $A$ equipped 
with an $n$-tuple 
$$\Delta_{\pi}^{(s)}=(\delta_{\pi,1}^{(s)},\ldots,\delta_{\pi,n}^{(s)})$$
 of higher $\pi$-derivations $A\rightarrow A$ of degree $s$. We will usually denote by
 $$\Phi_{\pi}^{(s)}=(\phi_{\pi,1}^{(s)},\ldots,\phi_{\pi,n}^{(s)})$$
 the attached family of
 higher $\pi$-Frobenius lifts.
 \end{definition}

\begin{remark}
In \cite{Bu95, Bu05, Bu17, BMP}   higher $\pi$-derivations  of degree $1$ were simply called {\bf $\pi$-derivations}. Higher $p$-Frobenius lifts of degree $1$ were simply called {\bf Frobenius lifts}. Finally partial $\delta_{\pi}$-rings of degree $1$ were called {\bf $\delta_{\pi}$-rings} (if $n=1$) or {\bf partial $\delta_{\pi}$-rings} (if $n\geq 2$).
\end{remark}

\begin{remark} For every $\pi\in \Pi$ the field $K_{\pi}$ is mapped into itself by 
every higher Frobenius automorphism $\phi^{(s)}$ of  $K^{\text{alg}}$ of degree $s\geq 1$ and the induced automorphism 
$\phi^{(s)}_{\pi}\colon R_{\pi}\ra R_{\pi}$   is a higher $\pi$-Frobenius lift of degree $s$ hence induces a higher $\pi$-derivation $\delta_{\pi}^{(s)}$ of degree $s$ on $R_{\pi}$. Therefore, if one is given a tuple $\Phi^{(s)}=(\phi^{(s)}_1,\ldots,\phi^{(s)}_n)$ of 
higher Frobenius automorphisms  of  $K^{\text{alg}}$ of degree $s\geq 1$ 
then for every $\pi\in \Pi$ the family  $\Phi^{(s)}_{\pi}=(\phi^{(s)}_{\pi,1},\ldots,\phi^{(s)}_{\pi,n})$ of restrictions of the members of $\Phi$ to $R_{\pi}$ defines a
structure of partial $\delta_{\pi}$-ring of degree $s$ on $R_{\pi}$  which we denote by 
$\Delta^{(s)}_{\pi}=(\delta^{(s)}_{\pi,1},\ldots,\delta^{(s)}_{\pi,n})$.
\end{remark}

 \begin{convention}\label{konvention1}
  From now on, until Section~\ref{sec:OC}, we will fix $\pi\in \Pi$.
To avoid a proliferation of subscripts $\pi$, we will  suppress $\pi$ in the notation throughout {\it as long as $\pi$ is fixed} and no confusion can arise. In particular, instead of 
$\delta^{(s)}_{\pi,i}, \phi_{\pi,i}^{(s)}, \Delta_{\pi}^{(s)}, \Phi_{\pi}^{(s)}$ we will write 
$\delta_i^{(s)},\phi_i^{(s)},\Delta^{(s)},\Phi^{(s)}$. Partial $\delta_{\pi}$-rings of degree $s$ will simply be referred to as {\bf partial $\delta$-rings} of degree $s$.
However we will not drop the index $\pi$ from $R_{\pi}$ because this would lead to ambiguities. Special attention to these subscripts will be paid when we considers multiple $\pi$'s at one time. This occurs when addressing overconvergence questions, see Section~\ref{sec:OC}. 
\end{convention}

\begin{convention}\label{konvention2}
 For the rest of this Section we fix, for each $s\geq 1$,  a structure of partial $\delta$-ring of degree $s$ on $R_{\pi}$ i.e.,  a family 
 $\Delta^{(s)}=(\delta_{1}^{(s)},\ldots,\delta_{n}^{(s)})$  of higher $\pi$-derivations of degree $s$, and we denote by 
$\Phi^{(s)}=(\phi^{(s)}_{1},\ldots,\phi^{(s)}_{n})$ the attached family of higher $\pi$-Frobenius lifts of degree $s$ on $R_{\pi}$.  Every structure of a partial $\delta$-ring  of degree $s$ on an $R_{\pi}$-algebra will be assumed compatible (in the obvious sense) with our fixed partial $\delta$-ring  structure of degree $s$ on $R_{\pi}$.
\end{convention}

\begin{notation} \label{MMn}
Let ${\mathbb M}_n$ be the free, non-commutative, monoid  with identity generated by
 the set $\{1,\ldots,n\}$. So the elements of $\mathbb M_n$ are the {\bf empty word} (denoted by $0$) together with {\bf words}  
$\mu:=i_1\ldots i_s$ where $i_1,\ldots,i_s\in \{1,\ldots,n\}$, $s\geq 1$.
The {\bf length}  of a word $\mu$ as above is defined by $|\mu|:=s$,
and the length of $0$ is defined by $|0|=0$. Multiplication is given by concatenation $(\mu,\nu)\mapsto \mu\nu$ and $0$ is the identity element. For all $r\in \mathbb N\cup\{0\}$ let ${\mathbb M}^r_n$ be the set of all elements in ${\mathbb M}_n$ of length $\leq r$. Set $\mathbb M_n^+:=\mathbb M_n\setminus \{0\}$, 
$\mathbb M_n^{r,+}:=\mathbb M_n^r\setminus \{0\}$, and $\mathbb M_n^{(r)}:=\mathbb M_n^r\setminus \mathbb M_n^{r-1}$.
For a  family of  elements $\phi_1,\ldots,\phi_n$ in a  monoid $\mathfrak G$ we have  a monoid homomorphism
$${\mathbb M}_n\rightarrow \mathfrak G,\ \ \mu=i_1\ldots i_l \mapsto \phi_{\mu}:=\phi_{i_1}\ldots \phi_{i_l},\ \ 0\mapsto 1.$$
\end{notation}

\begin{remark}\label{coburn}
We now specialize the above discussion to the case of the monoid $\mathfrak G$ of ring endomorphisms of a partial $\delta$-ring $A$ of degree $s$ with higher $\pi$-derivations  $\delta_1^{(s)},\ldots,\delta_n^{(s)}$ and associated higher $\pi$-Frobenius lifts $\phi_1^{(s)},\ldots,\phi_n^{(s)}\in \mathfrak G$. For all $\mu:=i_1\ldots i_l\in \mathbb M_n$
 we set 
$$(\phi^{(s)})_{\mu}:=\phi^{(s)}_{i_1} \circ \ldots 
\circ  \phi^{(s)}_{i_l}\in \mathfrak G   \textrm{ and }  (\delta^{(s)})_{\mu}:=\delta^{(s)}_{i_1} \circ \ldots \circ  \delta^{(s)}_{i_l}:A\rightarrow A.$$
Note that the maps $(\phi^{(s)})_{\mu}$ are higher $\pi$-Frobenius lifts of degree $ls$ and will be denoted simply by $\phi_{\mu}^{(ls)}$. However $(\delta^{(s)})_{\mu}$ are not higher $\pi$-derivations in general; rather, the  higher $\pi$-derivation associated to $\phi_{\mu}^{(ls)}$, which will be denoted by $\delta_{\mu}^{(ls)}$, can be expressed as a `polynomial' in the maps $(\delta^{(s)})_{\nu}$ with $|\nu|\leq |\mu|$.
For instance we have the following computation for $a\in A$:
$$\begin{array}{rcl}
\phi_{ij}^{(2s)}a & = &  \phi^{(s)}_{i}\phi^{(s)}_{j} a\\
\ & \ & \ \\
\ & = &    \phi^{(s)}_{i}(a^{p^s}+\pi\delta^{(s)}_{j}a)\\
\ & \ & \ \\
\ & = & (\phi^{(s)}_{i}(a))^{p^s}+\phi^{(s)}_{i}\pi\cdot \phi^{(s)}_{i}\delta^{(s)}_{j}a\\
 \ & \ & \ \\
\ & = & (a^{p^s}+\pi\delta^{(s)}_{i}a)^{p^s}+(\pi^{p^s}+\pi \delta^{(s)}_{i}\pi)((\delta^{(s)}_{j}a)^{p^s}+\pi\delta^{(s)}_{i}\delta^{(s)}_{j}a)\\
 \ & \ & \ \\
\ & \equiv & a^{p^{2s}}+\pi\cdot \delta^{(s)}_{i} \pi\cdot (\delta^{(s)}_{j}a)^{p^s}\ \ \textup{mod}\ \ \pi^2.\end{array}
$$
Hence, recalling that $(\delta^{(s)})_{ij}:=\delta^{(s)}_{i}\circ \delta^{(s)}_{j}$, we have:
\begin{equation}
\label{tasdi}
\delta_{ij}^{(2s)}a    =    \frac{1}{\pi} ((a^{p^s}+\pi\delta^{(s)}_{i}a)^{p^s}-a^{p^{2s}})
+(\pi^{p^s-1}+ \delta^{(s)}_{i}\pi)((\delta^{(s)}_{j}a)^{p^s}+\pi (\delta^{(s)})_{ij}a).
\end{equation}
\end{remark}

In what follows we record some basic formulae that play a role in the theory, and will play a particular role in Part 2 of this paper.

\begin{lemma}\label{uiui}
Let $A$ be a  partial $\delta$-ring of degree $s$ flat over $R_{\pi}$. For all $a\in A$ we have the following congruences in $A$:
$$\begin{array}{rcl}
\phi^{(s)}_{i}\phi^{(s)}_{j} a-\phi^{(s)}_{j}\phi^{(s)}_{i} a & \equiv  & 0\ \ \ \textup{mod}\ \ \pi,\\
\ & \ & \ \\
\delta_{ij}^{(2s)}a & \equiv  & \delta^{(s)}_{i} \pi\cdot (\delta^{(s)}_{j}a)^{p^s}
\ \ \textup{mod}\ \ \pi,\\
\ & \ & \ \\
\frac{1}{\pi}(\phi^{(s)}_{i}\phi^{(s)}_{j} a-\phi^{(s)}_{j}\phi^{(s)}_{i} a) & \equiv  & 
\delta^{(s)}_{i}\pi \cdot (\delta^{(s)}_{j} a)^{p^s}- \delta^{(s)}_{j}\pi \cdot (\delta^{(s)}_{i} a)^{p^s}\ \ \ \textup{mod}\ \ \pi\\
\ & \ & \ \\
\ & \equiv & \delta_{ij}^{(2s)}a-\delta^{(2s)}_{ji}a\ \ \ \textup{mod}\ \ \pi.\end{array}
$$
If in addition 
$$\phi^{(s)}_{i} a\equiv \phi^{(s)}_{j} a \ \ \ \textup{mod}\ \ \pi^2$$
for all $a\in R_{\pi}$, or equivalently  $\phi^{(s)}_i(\phi^{(s)}_j)^{-1}$ belongs to the first ramification group of the extension $K_{\pi}/K$, then for all $a\in R_{\pi}$ we have:
$$\delta^{(s)}_i a\equiv \delta^{(s)}_j a\ \ \textup{mod}\ \ \pi.$$
\end{lemma}

{\it Proof}. Follows immediately from the computation in Remark \ref{coburn}.
\qed

\medskip

We  record the structure of constants in our theory in the following lemma.

 \begin{lemma}
 \label{vinne}
 For  $a\in R_{\pi}$ the following hold:
 
 1) If $a$ is $0$ or a root of unity in $R$ then $\delta_i^{(s)}a=0$ for all $i$.
 
 2) If there exists $i$ such that $\delta_i^{(s)}a=0$ 
  then $a$ is either $0$ or a root of unity in $R$.
  
  3) If $a$ is a root of unity in $R$ and $b\in R_{\pi}$ then 
  $$\delta_i^{(s)}(ab)=a^{p^s}\delta_i^{(s)}b,\ \ \phi_i^{(s)}(ab)=a^{p^s}\phi_i^{(s)}b.$$
 \end{lemma}
 
 {\it Proof}. Assertions (1) and (3) follow directly from the definitions. For assertion (2)
 assume $\delta_i^{(s)}a=0$ (equivalently $\phi_i^{(s)}(a)=a^{p^s}$) and assume 
 $a$ is neither $0$ nor a root of unity in $R$; we will derive a contradiction. Since $\phi_i^{(s)}(a)$ has the same valuation as $a$ we have that $a\in R_{\pi}^{\times}$. Then we can write
 $a=\zeta b$ with $\zeta$ a root of unity in $R$ and $$1\neq b\in 1+\pi^{\nu}R_{\pi}$$ for some $\nu\in \mathbb N$; choose $\nu$ maximal  with this property.
 By assertion (1) we get $\phi_i^{(s)}(b)=b^{p^s}$. Writing $b=1+\pi^{\nu}c$  with $c\in R_{\pi}$ we get
 $$\phi_i^{(s)}(1+\pi^{\nu}c)=1+\phi_i^{(s)}(\pi)^{\nu}\phi_i^{(s)}(c)=(1+\pi^{\nu}c)^{p^s}\equiv 1\ \ \textup{mod}\ \ \pi^{\nu+1},$$
 which  implies $c\in \pi R_{\pi}$, contradicting the maximality of $\nu$.
\qed

\bigskip

For the next lemma recall that for every $a\in R_{\pi}$ there exist unique elements $\zeta_0,\zeta\in R$ each of which is  either $0$ or a root of unity in $R$  such that
\begin{equation}
\label{aeru}
a\equiv \zeta_0+\zeta\pi\ \ \ \textup{mod}\ \ \pi^2.
\end{equation}

\begin{lemma}\label{rou} With notation as in Equation \ref{aeru} the following formulae hold:
$$\begin{array}{rcl}
\delta_i^{(s)}a & \equiv &  \delta_i^{(s)}\pi \cdot \zeta^{p^s}\ \ \textup{mod}\ \ \pi\\
\ & \ & \ \\
\delta_{ij}^{(2s)} a & \equiv & \delta_i^{(s)} \pi \cdot (\delta_j ^{(s)}\pi)^{p^s}\cdot \zeta^{p^{2s}}\ \ \textup{mod}\ \ \pi.
\end{array}
$$
\end{lemma}

{\it Proof}.
The first congruence  follows by a trivial computation using the fact that
$\phi_i(\zeta_j)=\zeta_j^{p^s}$ for $j\in \{0,1\}$ and $i\in \{1,\ldots,n\}$.  The second congruence follows from the first congruence in our lemma plus the second congruence in Lemma \ref{uiui}.
\qed

\subsection{Complements on arithmetic ODEs} Recall that $R$ denotes the completed ring of integers of the maximal unramified extension of $\mathbb{Q}_p$ and has a unique Frobenius lift $\phi$. Set 
$$a^{\phi^s}:=\phi^s(a)=:a^{p^s}+p\delta^{(s)}a,\ \ \ a\in R.$$

\begin{lemma}\label{treiparti}
  For $a,b\in R$ with $a\equiv b$ mod $p^{\nu}$, $\nu\geq 1$, we have
$$a\equiv b\ \ \textup{mod}\ \ p^{\nu+1} \ \ \Leftrightarrow \ \ \delta^{(s)}a\equiv \delta^{(s)} b\ \ \textup{mod}\ \ p^{\nu}.$$ \end{lemma}

{\it Proof}. Write $a=b+p^{\nu}c$ with $c\in R$. We have
$$\begin{array}{rcl}
\delta^{(s)}a &\equiv  &\delta^{(s)} b+ \delta^{(s)}(p^{\nu}c) \ \ \textup{mod}\ \ p^{\nu}\\
\ & \equiv & \delta^{(s)} b+ p^{\nu-1}c^{p^s} \ \ \textup{mod}\ \ p^{\nu}
\end{array}$$
hence $c\in pR$ if and only if $\delta^{(s)}a\equiv \delta^{(s)} b$ mod $p^{\nu}$.\qed

\bigskip

We have the following existence and uniqueness result for first order ODEs.

\begin{proposition}\label{flowww}
Let 
$y,z$ be two $m$-tuples of indeterminates, consider restricted power series
$$F_1,\ldots,F_m,G\in R[y,z]^{\widehat{\ }},\ \ G\neq 0,$$ and set
$$f:=:F/G:=(F_1/G,\ldots,F_m/G).$$ Let $u^{(0)}\in R^m$ be such that 
$$G\left(u^{(0)},(u^{(0)})^{(p^s)}\right)\not\equiv 0\ \ \textup{mod}\ \ p.$$
 Then there exists a unique 
$u\in R^m$ such that $u\equiv u^{(0)}$ mod $p$ and 
\begin{equation}
\label{thesystemm}
\delta^{(s)}u= f\left(u,u^{\phi^s}\right).\end{equation}
\end{proposition}

{\it Proof}.
For the existence part it is enough to  construct a sequence $(u^{(\nu)})_{\nu\geq 1}$ of vectors $u^{(\nu)}\in R^n$ such that
$$\begin{array}{rcllll}
u^{(\nu+1)} & \equiv & u^{(\nu)} & \text{mod} & p^{\nu+1} &\textup{for}\ \ \nu\geq 0,\\
\ & \ & \ \\
\delta^{(s)}u^{(\nu)} & \equiv&  f(u^{(\nu)},(u^{(\nu)})^{\phi^s}) &  \textup{mod} & p^{\nu} & \textup{for}\ \ \nu\geq 0.\end{array}$$
Assume $u^{(\nu)}$ was constructed and set $u^{(\nu+1)}=u^{(\nu)}+p^{\nu+1}b$ with $b\in R^n$ to be determined. 
We have
$$f\left(u^{(\nu+1)},(u^{(\nu+1)})^{\phi^s}\right)\equiv f\left(u^{(\nu)},(u^{(\nu)})^{\phi^s}\right)\ \ \textup{mod}\ \ p^{\nu+1}.$$
On the other hand we have
$$\delta^{(s)}u^{(\nu)}=f\left(u^{(\nu)},(u^{(\nu)})^{\phi^s}\right)+p^{\nu}c$$
for some $c\in R^n$ so we get
$$\begin{array}{rcl}
\delta^{(s)}u_{\nu+1} & = & \delta^{(s)}u^{(\nu)}+\delta^{(s)}(p^{\nu+1}b)+C^{(s)}(u^{(\nu)},p^{\nu+1}b)\\
\ & \ & \ \\
\ & \equiv & f(u^{(\nu)},(u^{(\nu)})^{\phi^s})+p^{\nu}c+p^{\nu}b^{(p^s)}\ \ \textup{mod}\ \ p^{\nu+1}
\end{array}
$$
and it is enough to take $b$ such that $b^{(p^s)}\equiv c\ \ \textup{mod}\ \ p$ which is possible because the field $k$ is perfect.

For the uniqueness part 
if $u$ and $\tilde{u}$ are two vectors satisfying our equation  and congruent to $u^{(0)}$ mod $p$  then, using  Lemma \ref{treiparti}, one proves by induction on $\nu\geq 1$  that $u\equiv \tilde{u}$ mod $p^{\nu}$.\qed

\begin{corollary}
The map $\delta^{(s)}:pR\rightarrow R$ is bijective.\end{corollary}

 We denote the inverse of $\delta^{(s)}$ as $\int^{(s)}:R\rightarrow pR$. This can be viewed as {\bf $p$-integration operator}. It is, of course, non-additive but it is continuous for the $p$-adic topologies; in fact, by Lemma \ref{treiparti} if
   $a,b\in R$ and  $\nu\geq 1$ we have
$$a \equiv b\ \ \textup{mod} \ \ p^{\nu}\ \ \Leftrightarrow \int^{(s)} a\equiv \int^{(s)} b \ \ \textup{mod}\ \ p^{\nu+1}.$$

Finally we would like to `evaluate' solutions to ODEs at various points.
Note that  the scheme $\textup{Spec}(R)$ has a unique $k$-point over $R$, given by the canonical surjection
$$R\rightarrow k,\ \ a\mapsto a\ \ \textup{mod}\ \ p.$$
So if we need to talk about `other points' we need to generalize the notion of point
  and  consider the corresponding evaluation maps. We do this taking a clue from the Taylor expansion in classical calculus.
 
 \begin{definition}\label{evaeva}
 Let $t,t',t'',\ldots,t^{(i)},\ldots$ be indeterminates. The free $\delta$-ring $R\{t\}$ is the ring of polynomials
 $$R\{t\}:=R[t,t',t'',\ldots]$$
 and  will be called the ring of {\bf $\delta$-polynomials} in the variable $t$.
 For every $P\in R\{t\}$, $s\in \mathbb N$,
 and $a\in R$ we set
 $$P^{(s)}\{a\}:=P(a,\delta^{(s)} a,(\delta^{(s)})^2 a,\ldots)
 $$
 For $s$ fixed we 
 define the {\bf evaluation map at $P$},
 $$e_P^{(s)} \colon R\rightarrow k,\ \ a\mapsto a(P):=P^{(s)}\{a\}\ \ \textup{mod}\ \ p.$$ When the degree $s$ is clear, we simply write this as $e_P$.
 Set $P_0 :=t$ and note that $e_{P_0}$ is the canonical surjection $R\rightarrow k$ and we have the formula
 \begin{equation}
 \label{Pada}
 a(P)=(P^{(s)}\{a\})(P_0),\ \ a\in R.\end{equation}
  \end{definition}
  
  \begin{remark}\ 
  
  1) 
  We  view $\delta$-polynomials as `generalized points' of $R$. 
  
  2) 
  If $P, Q \in R\{t\}$ 
  satisfy $P\equiv Q$ mod $p$ then clearly $e_P=e_Q$. Conversely, if $s=1$ and 
  $P,Q \in R\{t\}$ satisfy $e_P=e_Q$ then $P\equiv Q$ mod $p$; this follows from 
  \cite[Lemma 3.20]{Bu05}.
  
  3) If $P^{(s)}\{a\}=0$ for all $a\in R$ then $P=0$; the argument when $s=1$ follows from  \cite[Lem. 3.20]{Bu05} and the general case follows easily from the $s=1$ case.
  
  4) 
  The evaluation of $a$ at  $P$ is defined to be the evaluation at $P_0$
  of
  $P^{(s)}\{a\}$ which is a situation that mimics that in the case of Taylor expansions, where, at least formally,  the role of $P$ is played by a `linear differential operator of infinite order;' see Remark \ref{uff}.
  In our case here the order of every  $P$ is, of course, finite. 
  
  5) For $a=(a_1,\ldots,a_n)\in R^n$ and $P\in R\{t\}$ one can define
  $$a(P):=(a_1(P),\ldots,a_n(P))\in k^n.$$
  
  6) 
One can view the ring $R\{t\}$ a a graded ring by giving $t^{(i)}$ the degree
$$\deg(t^{(i)}):=p^{is}.$$ 
One can easily check, using Lemma \ref{vinne}, assertion (3), that a  $\delta$-polynomial $P\in R\{t\}$ is homogeneous of degree $d$ if and only if for all $a\in R$ and every root of unity $\zeta \in R$ we have
$$P^{(s)}\{\zeta a\}=\zeta^d P^{(s)}\{a\}.$$
\end{remark}
  
 \begin{remark}\label{uff}The $\delta$-polynomial $P_0=t$ can be viewed as  an analogue of the `origin' $0$ in an open interval $I\subset \mathbb R$ containing $0$.
  One is tempted to single out a class of $\delta$-polynomials $P\in R\{t\}$ that are analogues of arbitrary points in $I$. Here is a candidate for such a class of $\delta$-polynomials that we shall call {\bf Taylor $\delta$-polynomials}. Let $s=1$ and let us suppress the superscript $(s)$.
  By \cite[Prop. 2.10]{BPS} that there exist $\delta$-polynomials 
  $P_0,P_1,P_2,\ldots$ where $$P_i\in \mathbb Z_p[t,t',\ldots,t^{(i)}],\ \ i\geq 0$$
  with $P_0=t$ such that
  the natural map
  \begin{equation}
  \label{natutu}
  R\rightarrow W(R)=R^{\mathbb N}\end{equation}
  to the ring of Witt vectors $W(R)$ induced by the  Frobenius lift on $R$  is given by
  $$a\mapsto (P_0\{a\},P_1\{a\},P_2\{a\},\ldots), \ \ a\in R.$$
    There is an explicit recipe to construct these $\delta$-polynomials; in particular
    $P_i$ is homogeneous of degree $p^i$ and we have
    $$P_i-t^{(i)}\in \mathbb Z_p[t,t',\ldots,t^{(i-1)}],\ \textup{for}\ \  i\geq 1.$$
    A $\delta$-polynomial $P\in R\{t\}$ will be called a {\bf Taylor $\delta$-polynomial} if it has the form
  $$P:=\sum_{i=0}^{\infty} a_i P_i\in R\{t\},$$
  with $a_i\in R$  almost all $0$. If more that one $a_i$ is non-zero such a polynomial is, of course, not homogeneous.
  
  The terminology above is justified by the following analogy with calculus.
    Indeed the map (\ref{natutu}) is an analogue of the `exponential' operator
  $$L\rightarrow L[[t]],\ \ \lambda\mapsto e^{tD}\lambda:=
  \sum_{i=0}^{\infty}\frac{t^iD^i \lambda}{i!},$$
  for a field of characteristic zero $L$ equipped with a derivation $D:L\rightarrow L$; the partial sums in $e^{tD}\lambda$ are the usual Taylor polynomials in calculus.
  
  Using \cite[Lem. 3,20]{Bu05} one trivially has that if $a,b\in R$ are such that $a(P)=b(P)$ for all Taylor $\delta$-polynomials $P\in R\{t\}$, intuitively if the values of $a$ and $b$  at all Taylor $\delta$-polynomials  coincide, then $a=b$. \end{remark}

 As an application of the above constructions
 assume $f$ is as  in Proposition \ref{flowww}. Fix  $P\in R\{t\}$, let $S\subset k^m$ be the zero set of the polynomial $G(y,y^{(p^s)})$,  and consider the {\bf transport} map 
 $$\textup{trans}_{f,P}:k^m\setminus S\rightarrow k^m$$ defined as follows: for every $\lambda^{(0)}\in k^m\setminus S$ let $u\in R^m$ be the unique vector such that
 $$\begin{array}{rcl}\delta^{(s)}u &= &f(u,u^{\phi^s}),\\
 \  \ & \ \\
 u(P_0) & = & \lambda^{(0)};\end{array}$$
 cf. Proposition \ref{flowww}.
 Then define
 $$\textup{trans}_{f,P}(\lambda^{(0)}):=u(P).$$
 Let $g\in k[y]$ be the image of $G(y,y^{(p)})$ and consider the schemes
 $$\mathbb A_k^m:=\textup{Spec}(k[y]),\ \ \mathbb S:=\textup{Spec}(k[y]/(g)).$$
 Furthermore view $k^m$ and $k^m\setminus S$ as the sets of $k$-points of the schemes
 $\mathbb A_k^m$ and $\mathbb A_k^m\setminus \mathbb S$, respectively.
 
 \begin{proposition}\label{regularr}
 The map $\textup{trans}_{f,P}$ is induced by a morphism of $k$-schemes (still denoted by)
 $$\textup{trans}_{f,P}:\mathbb A_k^m\setminus \mathbb S\rightarrow \mathbb A_k^m.$$
 \end{proposition}
 
 {\it Proof}. Set $f^{(1)}:=f$. Applying $\delta^{(s)}$ to the equality 
 \begin{equation}
 \label{skippp}
 \delta^{(s)}u =f^{(1)}(u,u^{\phi^s})\end{equation}
 one gets an equality of the form
 \begin{equation}
 \label{skipppp}
 (\delta^{(s)})^2u =f^{(12)}(u,\delta^{(s)}u,(\delta^{(s)}u)^{\phi^s}),\end{equation}
 where $f^{(12)}=F^{(12)}/G^{(12)}$, $F^{(12)}$ is an $m$-tuple of polynomials in $3m$ variables with $R$-coefficients, and
 $G^{(12)}$ is a polynomial in $3m$ variables with $R$-coefficients with $G^{(12)}$ congruent mod $p$ to a power of $G$. Substituting (\ref{skippp}) into (\ref{skipppp}) we get an equality
 \begin{equation}
 \label{skippppu}
 (\delta^{(s)})^2u =f^{(2)}(u,u^{\phi^s},u^{\phi^{2s}}),\end{equation}
 where $f^{(2)}=F^{(2)}/G^{(2)}$, $F^{(2)}$ is an $m$-tuple of polynomials in $3m$ variables with $R$-coefficients, and
 $G^{(2)}$ is a polynomial in $3m$ variables with $R$-coefficients with $G^{(2)}$ congruent mod $p$ to a power of $G$. Repeating this procedure we get, for all $i\geq 2$ an equality of the form
  \begin{equation}
 \label{skippppuw}
 (\delta^{(s)})^iu =f^{(i)}(u,u^{\phi^s},\ldots,u^{\phi^{is}}),\end{equation}
 where $f^{(i)}=F^{(i)}/G^{(i)}$, $F^{(i)}$ is an $m$-tuple of polynomials in $(i+1)m$ variables with $R$-coefficients, and
 $G^{(i)}$ is a polynomial in $(i+1)m$ variables with $R$-coefficients with $G^{(i)}$ congruent mod $p$ to a power of $G$. If $f^{(i)}$ has components $f^{(i)}_k$
 we get
 $$(\delta^{(s)})^i u_k \equiv f_k^{(i)}(u,u^{(p^s)},\ldots,u^{(p^{is})})\ \ \ \textup{mod}\ \ p.$$
 Hence
 $$u_k(P)=P(u,f_k^{(1)}(u,u^{(p^s)}),f_k^{(2)}(u,u^{(p^s)},u^{(p^{2s})}),\ldots) \ \ \textup{mod}\ \ p.$$
 \qed

\section{Arithmetic connections in the PDE setting}
   
   An arithmetic ODE Levi-Civita formalism was developed in \cite[Sec. 4.3]{Bu17}.
   An extension of that theory to several arithmetic directions was developed in \cite{Bu19}, however this was done by introducing new metrics and  was reducible to an ODE theory. We now address  genuinely PDE versions of the Levi-Civita and Chern connections.
   
    Throughout this section we adopt the Conventions \ref{konvention1} and \ref{konvention2}. In particular, we recall that we have fixed $\pi\in \Pi$ and for every $s\geq 1$ we have fixed a partial $\delta$-ring structure of degree $s$ on $R_{\pi}$; all partial $\delta$-ring structures of degree $s$ will be assumed  compatible with the above fixed one on $R_{\pi}$. The numbers (\ref{4numbers}), $e,s,n,N$, are not assumed to be related.

   \subsection{$\pi$-connections and Christoffel symbols}
   We start by recalling/adapting some basic notions from \cite{Bu17, Bu19}.

  \begin{notation} We begin with some matrix notation and conventions that will be repeatedly used later.
For every ring $S$ and every integer $N\geq 1$ we denote by $\textup{Mat}_N(S)$ the ring of $N\times N$ matrices with coefficients in $S$ and by  $\textup{GL}_N(S)$ the group of invertible elements of $\textup{Mat}_N(S)$.
 We denote by $\textup{Mat}_N(S)^{\textup{sym}}$ and $\textup{GL}_N(S)^{\textup{sym}}$ the sets of symmetric matrices in $\textup{Mat}_N(S)$ and $\textup{GL}_N(S)$, respectively. The identity matrix  will  be denoted by $1=1_N=(\delta_{kl})$ where $\delta_{kl}$ is the Kronecker symbol. Using $\delta$ to denote both the Kronecker symbol and our $\pi$-derivations should not lead to any confusion: the meaning should be clear each time from context or recalled as needed.

 For every  $X\in \textup{Mat}_N(S)$  we denote by $X_{ij}$ its entries and we write $X=(X_{ij})$. Recall that for  every integer $d\geq 1$ we denote by $X^{(p^d)}$ the matrix $(X_{ij}^{p^d})$ and by $X^t$ the transpose of $X$. Similarly for $v=(v_i)$  a row or column vector with entries $v_i\in S$ we denote by $v^{(p^d)}$ the vector with entries $v_i^{p^d}$.  If $u:S\rightarrow S'$ is a map of sets we write $u(X):=(u(X_{ij}))$ and $u(v)=(u(v_i))$. In what follows, for $N=n$, we will be led to consider families 
 $(X_1,\ldots,X_n)$ of matrices in $\textup{Mat}_n(S)$. For such a family we denote by $(X_i)_{jk}$ the  $jk$-entries of $X_i$ so we write $X_i=((X_i)_{jk})$. However we will {\it never} use the notation $X_{ijk}$ to denote the $jk$-entry of $X_i$. This notation  will be reserved for objects obtained by  `lowering the indices,' in a sense that will be explained  later.  On the other hand we will also consider matrices  $X^k$ indexed by superscripts $k$, rather than subscripts. As a rule these superscripts will never mean `raising $X$ to the $k$-th power.'
 For a matrix $X^k$ we denote by $X^k_{ij}:=(X^k)_{ij}$ its $ij$-entries. For every $n$-tuple of matrices $(X_1,\ldots,X_n)$ we  define the $n$-tuple of matrices $(X^1,\ldots,X^n)$ by the formula 
\begin{equation}\label{treimat}
(X^k)_{ij}=X^k_{ij}:=(X_i)_{jk}=(X_i^t)_{kj};\end{equation}
 and, vice versa, if one is given an $n$-tuple of matrices $(X^1,\ldots,X^n)$ the above formula defines an $n$-tuple of matrices
 $(X_1,\ldots,X_n)$.
\end{notation}
       
   \begin{notation}\label{notynoty} Let $x=(x_{ij})$ be an $N\times N$  matrix of indeterminates. We consider the group scheme
 $$G=\textup{GL}_{N,R_{\pi}}=\textup{Spec }R_{\pi}[x,\det(x)^{-1}],$$
 so for the ring of global functions we have 
$\mathcal O(G)=R_{\pi}[x,\det(x)^{-1}]$ and for the group of points we have $G(R_{\pi})=\textup{GL}_N(R_{\pi})$. 
We consider the ring
$$\mathcal A:=\widehat{\mathcal O(G)}=R_{\pi}[x,\det(x)^{-1}]^{\widehat{\ }}$$
and we consider the prime ideal of $\mathcal A$,
$$\mathcal P:=(x-1),$$
generated by  the entries $x_{kl}-\delta_{kl}$ of the matrix $x-1$. We also consider the 
maximal ideal of $\mathcal A$,
$$\mathcal M:=(\pi,\mathcal P)=(\pi,x-1),$$
 generated by $\pi$ and $\mathcal P$.
 The residue field $\mathcal A/\mathcal M$ equals $k$. Also we have a canonical ring homomorphism given by the substitution $x\mapsto 1$,
$$\mathcal A\rightarrow \mathcal A/\mathcal P\simeq R_{\pi},\ \ \ F\mapsto F(1):=F_{|x=1}.$$
\end{notation}
 
 \begin{definition}\label{defconnn}
A {\bf $\pi$-connection} of degree $s$ on $G$  is an $n$-tuple 
 $$\Delta^{(s)G}=((\delta^{(s)}_1)^G,\ldots,(\delta^{(s)}_n)^G)$$ 
 of higher $\pi$-derivations of degree $s$ on $\mathcal A$ extending the given higher $\pi$-derivations of degree $s$ on $R_{\pi}$, respectively. To give a 
$\pi$-connection  of degree $s$  on $G$ is the same as to give a structure of  partial $\delta$-ring of degree $s$ 
on 
$\mathcal A$.
We denote by $$\Phi^{(s)G}=((\phi^{(s)}_1)^G,\ldots,(\phi^{(s)}_n)^G)$$ the attached  family of higher 
$\pi$-Frobenius lifts of degree $s$.
\end{definition}

A `trivial' example is given as follows.

\begin{definition}\label{tritri}
The {\bf  trivial $\pi$-connection} of degree $s$, 
$$\Delta_0^{(s)G}=((\delta_{1,0}^{(s)})^G,\ldots,(\delta_{n,0}^{(s)})^G),$$  
is defined by $(\delta_{i,0}^{(s)})^G x_{jk}=0$.
It has attached higher $\pi$-Frobenius lift $(\phi_{i,0}^{(s)})^{G}$ satisfying $(\phi_{i,0}^{(s)})^G (x)=x^{(p^s)}$. \end{definition}

\begin{definition}\label{christos}
 The {\bf Christoffel symbol of the second kind} of a 
 $\pi$-connection  $((\delta_1^{(s)})^G,\ldots,(\delta_n^{(s)})^G)$  of degree $s$ on $G$ 
  is the $n$-tuple of matrices $$\Gamma^{(s)}=
  (\Gamma_1^{(s)},\ldots,\Gamma_n^{(s)})$$ whose transposed $\Gamma_i^{(s)t}$ are given by
  \begin{equation}
\label{penu}\Gamma_i^{(s)t}:=(x^{(p^s)})^{-1} (\delta_i^{(s)})^G x \in \textup{Mat}_N(\mathcal A).\end{equation}
In view of  Equation (\ref{treimat})  we write 
\begin{equation}
\label{trematt}
\Gamma_{ij}^{k(s)}:=(\Gamma_i^{(s)})_{jk}=(\Gamma_i^{(s)t})_{kj}\in \mathcal A.\end{equation}
Note that we will {\it not} use the notation $\Gamma^{(s)}_{ijk}$ to denote the $jk$ entry of $\Gamma^{(s)}_i$. The notation $\Gamma^{(s)}_{ijk}$ will be reserved for the `Christoffel symbols of the first kind' to be introduced later via `lowering of indices.' We have the following formula
\begin{equation}
\label{la11}
((\phi_i^{(s)})^G(x_{kj}))(1)=\delta_{kj}+\pi\Gamma^{k(s)}_{ij}(1),
\end{equation}
where $\delta_{kj}$ is the Kronecker symbol.

Finally we set $\Lambda^{(s)}:=(\Lambda_1^{(s)},\ldots,\Lambda_n^{(s)})$,
where
\begin{equation}
\Lambda_i^{(s)}:=(x^{(p^s)})^{-1}(\phi_i^{(s)})^G (x)=1+\pi \Gamma_i^{(s)t}
\in \textup{GL}_N(\mathcal A).\end{equation}
\end{definition}

For an `intrinsic' description of Christoffel symbols we refer to Subsection \ref{ld}.
For a discussion of the analogy with classical differential geometry we refer to  the Appendix of Part 2, \cite{BMadg2}. 

\begin{remark}
Our Christoffel symbols here are slightly different from the ones in \cite{Bu17} and \cite{Bu19};
the change was necessary in order to accommodate the new PDE setting. \end{remark}

\subsection{Parallel transport and geodesics}
Assume in this subsection only that the numbers (\ref{4numbers}), $e,s,n,N$ satisfy $e=n=N$.

\begin{definition}
The elements  of the $R$-module $\textup{Hom}_{R-\textup{mod}}(R_{\pi},R)$ will be referred to as {\bf curves} in $R_{\pi}$. The elements  of the $k$-linear space $\textup{Hom}_{R-\textup{mod}}(R_{\pi},k)$ will be referred to as {\bf points} of $R_{\pi}$. \end{definition}

For a discussion of the analogy with classical differential geometry see the Appendix of Part 2, \cite{BMadg2}.
 We fix an $R$-module basis $\theta:=(\theta_1,\ldots,\theta_n)$ of $R_{\pi}$. There is an identification 
\begin{equation}
\label{curcur}
R^n\simeq \textup{Mat}_{n\times 1}(R)\simeq \textup{Hom}_{R-\textup{mod}}(R_{\pi},R),\ \ \ c\mapsto 
c^*,\end{equation}
$$
c^*\left(\sum_{i=1}^n a_i\theta_i\right):=\sum_{i=1}^n a_ic_i,\ 
\ \ \textup{for}\ \ c:=(c_1,\ldots,c_n)^t,\ \ a_i\in R.$$
So we identify  curves $c^*$ with  column vectors $c\in R^n$. 
Intuitively, for $\theta\in R_{\pi}$ we view $c^*(\theta)$ as the restriction of the `function' $\theta$ to the curve $c$.
Also there is an identification
\begin{equation}
\label{curcur2}
k^n\simeq \textup{Mat}_{n\times 1}(k)\simeq \textup{Hom}_{R-\textup{mod}}(R_{\pi},k),\ \ \ \lambda\mapsto 
\lambda^*,\end{equation}
$$
\lambda^*\left(\sum_{i=1}^n a_i\theta_i\right):=\sum_{i=1}^n a_i\lambda_i,\ 
\ \ \textup{for}\ \ \lambda:=(\lambda_1,\ldots,\lambda_n)^t,\ \ a_i\in R.$$
So we identify  points $\lambda^*$ with  column vectors $\lambda\in k^n$. Intuitively, for $\theta\in R_{\pi}$ we view  
$$\theta(\lambda):=\lambda^*(\theta)$$
 as the value of the `function' $\theta$ at the point $\lambda$.

Recall, cf. Definition \ref{evaeva}, that we denoted by $e_{P_0}:R\rightarrow k$, $a\mapsto a(P_0)$ the reduction modulo $p$ map, viewed as the evaluation map at $P_0=t$ in $R\{t\}$.
This map induces a map
$$e_{P_0}:R^n\rightarrow k^n,\ \ \ c=(c_1,\ldots,c_n)^t\mapsto c(P_0):=(c_1(P_0),\ldots,c_n(P_0))^t.$$
Intuitively we view $c(P_0)$ as the  `origin' of the curve $c$.
We have the following compatibility between the above evaluation maps: 
\begin{equation}
\label{lowonbread}
\theta(\lambda)=(c^*(\theta))(P_0)\ \ \ \textup{for}\ \ c\in R^n,\ \ \lambda=c(P_0),\ 
 \theta\in R_{\pi}.
\end{equation}
Intuitively, if the origin of the curve $c$ is $\lambda$ then the `function' $\theta$ evaluated at $\lambda$ equals the restriction of $\theta$ to $c$  evaluated at $P_0$.

We write
$$\frac{p}{\pi}=\sum_{i=1}^nr_i\theta_i,\ \ r_i\in R.$$
Since $\frac{p}{\pi}\not\in pR_{\pi}$ we have
\begin{equation}
\label{dsa}
(r_1,\ldots,r_n)\not\in pR^n.
\end{equation} 

Recall that we denote by $\phi$  the Frobenius lift on $R$ and we write
$$a^{\phi^s}=\phi^s(a)=a^{p^s}+p\delta^{(s)}a$$ for $a\in R$, $\mathcal A = R_\pi[x, \det (x)^{-1}]$, and $\mathcal P = (x-1) \subset \mathcal A$.  For a curve 
$c\in R^n$  the composition 
\begin{equation}
\label{formerrho}
\mathcal A\rightarrow \mathcal A/\mathcal P\simeq R_{\pi}
\stackrel{c^*}{\longrightarrow} R,\ \ \ F\mapsto F_{c}:=c^*(F(1))\end{equation}
can be intuitively viewed as the `restriction of $F$ to $c$.'
Also, for $c\in R^n$ we consider the column vector
\begin{equation}
v=(v_1,\ldots,v_n)^t:=\delta^{(s)}c:=(\delta^{(s)}c_1,\ldots,\delta^{(s)}c_n)^t\in R^n\end{equation}
which we refer to   as the {\bf velocity} of $c$ (or of $c^*$).
For two curves $c,c^{(0)}\in R^n$ with $c(P_0)=c^{(0)}(P_0)$ and velocities $v,v^{(0)}$ we will  later consider the condition $v(P_0)=v^{(0)}(P_0)$, equivalently,
$$\delta^{(s)}c\equiv \delta^{(s)}c^{(0)}\ \ \ \textup{mod}\ \ p;$$
the latter condition can be thought of as $c$ and $c^{(0)}$ being 
`tangent at their common origin.'

Consider  a $\pi$-connection $\Delta^{(s)G}=((\delta_1^{(s)})^G,\ldots,(\delta_n^{(s)})^G)$ of degree $s$ on $G$
and a curve $c\in R^n$ with velocity $v$. We will consider the map
$$\delta_{v}^{(s)G}:=\sum_{i=1}^n v_i^{\phi^s}(\delta_i^{(s)})^G:\mathcal A\rightarrow \mathcal A.$$

\begin{definition}
Given a curve $c\in R^n$ with velocity $v\in R^n$ and a vector $w\in R^n$ 
we define the {\bf derivative} of $w$ along $c$ 
with respect to the $\pi$-connection $\Delta^{(s)G}$
to be the vector
$$w'_{c}:=(\delta_{v}^{(s)G}(xw))_{c}.$$
We say that $w$ is {\bf parallel along $c$} with respect to the $\pi$-connection $\Delta^{(s)G}$ if
\begin{equation}
\label{parabubu}
w'_{c}=0.
\end{equation}
\end{definition}

In the above formulae $w$ is viewed, as usual,  as a column vector with entries in $R$, hence
$xw$ is  a column vector with entries in $\mathcal A$,
and we view
$\delta_{v}^{(s)G}(xw)$ as a column vector with entries in $\mathcal A$ whose entries are obtained by applying $\delta_{v}^{(s)G}$ to the entries of $xw$. Thus 
$(\delta_{v}^{(s)G}(xw))_{c}$
is a column vector with entries in $R$ whose entries are obtained by applying the map (\ref{formerrho}) to the entries of $\delta_{v}^{(s)G}(xw)$.

The pair $(c,w)$ can be viewed as an analogue of a vector field along $c$.
For an arbitrary curve  $c\in R^n$ with velocity $v$ the vector
$$v'_{c}=(\delta_{v}^{(s)G}(xv))_{c}\in R^n$$  can be viewed as an analogue of `acceleration' of (or along) $c$.

\begin{remark}
The definition of `acceleration' $v'_{c}$ has the following invariance with respect to the action of the symmetric group $\Sigma_n$, which morally plays the role of the group of `coordinate changes' in our theory. In order to indicate the dependence on the bases 
we write $c^*_{\theta}$ instead of $c^*$ and also $(\delta_{v}^{(s)G}(xv))_{c,\theta}$ instead of $(\delta_{v}^{(s)G}(xv))_{c}$.
For every $\epsilon\in \Sigma_n$  let $P_{\epsilon}$ be the permutation matrix corresponding to $\epsilon$, that is $(P_{\epsilon})_{ij}=\delta_{\epsilon(i)j}$. Let $\tilde{x}:=xP^t_{\epsilon}$, $\tilde{\theta}:=P_{\epsilon}\theta$, and 
$\tilde{\Delta}^{(s)G}=\Delta^{(s)}P_{\epsilon}^t$; so the $i$-th component $(\tilde{\delta}_i^{(s)})^G$ of $\tilde{\Delta}^{(s)G}$ equals $(\delta_{\epsilon(i)}^{(s)})^G$. Let $c,\tilde{c}\in R^n$ define the same curve with respect to the bases $\theta$ and $\tilde{\theta}$, respectively; i.e., $c^*_{\theta}=\tilde{c}^*_{\tilde{\theta}}$. By our definition we have $c^*_{\theta}(\theta)=c$ and 
$\tilde{c}^*_{\tilde{\theta}}(\tilde{\theta})=\tilde{c}$. We deduce that $\tilde{c}=P_{\epsilon}c$, hence the velocities $v$ and $\tilde{v}$ of $c$ and $\tilde{c}$, respectively,  are related by the equation $\tilde{v}=P_{\epsilon}v$. On the other hand since $P_{\epsilon}P_{\epsilon}^t=1$ we have 
$$\delta_{\tilde{v}}^{(s)G}=(\tilde{v}^{\phi^s})^t \tilde{\Delta}^{(s)Gt}=(v^{\phi^s})^tP^t_{\epsilon}
\tilde{\Delta}^{(s)Gt}=(v^{\phi^s})^t\Delta^{(s)Gt}=\delta^{(s)G}_v.$$
Since $xv=\tilde{x}\tilde{v}$ we deduce an equality of the corresponding accelerations:
$$(\delta_{\tilde{v}}^{(s)G}(\tilde{x}\tilde{v}))_{\tilde{c},\tilde{\theta}}=(\delta_{v}^{(s)G}(xv))_{c,\theta}.
$$
\end{remark}

\begin{definition}
We say that the curve $c\in R^n$ with velocity $v$ is a {\bf geodesic} for the $\pi$-connection $\Delta^{(s)G}$ if 
\begin{equation}
\label{geo}
v'_{c}=0;
\end{equation}
in other words if $v$ is parallel along $c$ or rather that the acceleration vanishes.
\end{definition}

Recall that we denoted by  $v_i$ the components of $v$ and let $w_i$ be the components of $w\in R^n$.
Explicitly the $k$-th component 
$(w'_{c})_k$ of $w'_{c}$ is given by the following calculation:
$$
\begin{array}{rcl}
(w'_{c})_k & = & c^*\left(\left(
\sum_i v_i^{\phi^s}(\delta_i^{(s)})^G\left(\sum_j x_{kj}w_j\right)\right)(1)\right)\\
\ & \ & \ \\
\ & = &  c^*\left(\frac{1}{\pi}\left(
\sum_i v_i^{\phi^s}\left((\phi_i^{(s)})^G\left(\sum_j x_{kj}w_j\right)-\left(\sum_j x_{kj}w_j
\right)^{p^s}\right)\right)(1)\right)\\
\ & \ & \ \\
\ & = & c^*\left(\frac{1}{\pi}\left(
\sum_{ij} v_i^{\phi^s}w_j^{\phi^s}\left(\delta_{kj}+\pi\Gamma_{ij}^{k(s)}(1)\right)-
\sum_{ij}v_i^{\phi^s}w_j^{p^s}\delta_{kj}\right)\right),\ \ \textup{by (\ref{la11})}
\\
\ & \ & \ \\
\ & = & c^*\left(\frac{1}{\pi}\left(
\left(\sum_{i} v_i^{\phi^s}\right)w_k^{\phi^s}+\pi\sum_{ij}v_i^{\phi^s}w_j^{\phi^s}\Gamma_{ij}^{k(s)}(1)
-\left(\sum_{i} v_i^{\phi^s}\right)w_k^{p^s}
\right)\right)\\
\ & \ & \ \\
\ & = &  c^*\left(\frac{p}{\pi}
\left(\sum_{i} v_i^{\phi^s}\right)\delta^{(s)}w_k+\sum_{ij}v_i^{\phi^s}w_j^{\phi^s}\Gamma_{ij}^{k(s)}(1)
\right)\\
\ & \ & \ \\
\ & = &  c^*\left(\frac{p}{\pi}\right)
\left(\sum_{i} v_i^{\phi^s}\right)\delta^{(s)}w_k+\sum_{ij}c^*\left(\Gamma_{ij}^{k(s)}(1)\right)v_i^{\phi^s}w_j^{\phi^s}.
\end{array}
$$

\begin{definition} \label{pppo} Recall $p/\pi = \sum r_\ell \theta_\ell$. A curve $c\in R^n$  with velocity $v$
is  {\bf non-degenerate} if we have
$$\sum_l r_lc_l\not\equiv 0\ \ \ \textup{mod}\ \ \ p,$$
$$\sum_l v_l\not\equiv 0\ \ \ \textup{mod}\ \ \ p.$$ The first condition guarantees $c^*\left(\frac{p}{\pi}\right)$ is not zero modulo $p$.
\end{definition}

\begin{remark} By (\ref{dsa}) the set of non-degenerate curves in $R^n$ is non-empty; in fact a `generic' curve (in a sense that can be easily made precise) is non-degenerate. \end{remark}

The  computation preceding Definition \ref{pppo} thus shows the following.

\begin{proposition}\label{eqforgeo}
Let $c\in R^n$ be a  non-degenerate curve with velocity $v$ and write 
$$\Gamma^{k(s)}_{ij}(1)=\sum_{l=1}^n \Gamma^{k(s)}_{ij,l}\theta_l,\ \ 
\Gamma^{k(s)}_{ij,l}\in R.$$
Then a vector $w\in R^n$ is parallel along $c$ if and only if
\begin{equation}\label{pararara}
\delta^{(s)}w_k+\frac{\sum_{ij}(\sum_l \Gamma^{k(s)}_{ij,l}c_l)
v_i^{\phi^s}w_j^{\phi^s}}{(\sum_l r_l c_l)(\sum_l v_l^{\phi^s})}
=0,\ \ k\in \{1,\ldots,n\}.
\end{equation}
In particular  $c$ is a geodesic if and only if
 \begin{equation}\label{geeoo}
\delta^{(s)}v_k+\frac{\sum_{ij}(\sum_l \Gamma^{k(s)}_{ij,l}c_l)
v_i^{\phi^s}v_j^{\phi^s}}{(\sum_l r_l c_l)(\sum_l v_l^{\phi^s})}
=0,\ \ k\in \{1,\ldots,n\}.
\end{equation}
\end{proposition}

\begin{remark}\label{syyster}\ 

1) The system (\ref{pararara}) can be viewed as a `first order' system in the unknowns $w_1,\ldots,w_n$
and has the form (\ref{thesystemm}) considered in Proposition~\ref{flowww}. Note that if $w\in R^n$ is a solution to the system (\ref{pararara}) and $\zeta$ is a root of unity in $R$ then $\zeta w\in R^n$ is also a solution to the system (\ref{pararara}). This follows directly from Lemma \ref{vinne}, assertion (3). 
Note also that the system (\ref{pararara}) can be written in the form 
\begin{equation}
\label{fformma}
\delta^{(s)} w=\alpha w^{\phi^s},\end{equation}
for some matrix $\alpha\in \textup{Mat}_n(R)$ depending on $c$. In its turn the system (\ref{fformma}) can be rewritten in the form
\begin{equation}
\label{fformm}
\delta^{(s)} w=\beta w^{(p^s)},\end{equation}
where
$$ \beta:=\alpha+p\alpha^2+p^2\alpha^3+\ldots$$
For systems 
of the form (\ref{fformm}) with $s=1$ and arbitrary $\beta\in \textup{Mat}_n(R)$
a Galois theory was developed in \cite[Sec. 5.4]{Bu17}. It would be interesting to
investigate the Galois theoretic properties of the system (\ref{fformm}) arising from 
the system (\ref{pararara}).

2) The system consisting of the equations (\ref{geeoo}) together with the equations
\begin{equation}\label{suppp}
\delta^{(s)}c_k=v_k,\ \ k\in \{1,\ldots,n\}
\end{equation}
can be viewed as a system in the unknowns $c_1,\ldots,c_n,v_1,\ldots,v_n$ and has, again, the form (\ref{thesystemm}); this system is equivalent to 
the  `second order' system in the unknowns $c_1,\ldots,c_n$ given by:
 \begin{equation}\label{eqforgeooo}
(\delta^{(s)})^2(c_k)+\frac{\sum_{ij}(\sum_l \Gamma^{k(s)}_{ij,l}c_l)
(\delta^{(s)}c_i)^{\phi^s}(\delta^{(s)}c_j)^{\phi^s}}{(\sum_l  r_lc_l)(\sum_l (\delta^{(s)}c_l)^{\phi^s})}
=0,\ \ k\in \{1,\ldots,n\}.
\end{equation}
 Note that if $c\in R^n$ is a solution to the system (\ref{eqforgeooo}) and $\zeta$ is a root of unity in $R$ then $\zeta c\in R^n$ is also a solution to the system (\ref{eqforgeooo}).
 Again, this follows directly from Lemma \ref{vinne}, assertion (3).
\end{remark}

\begin{remark}
Assume 
\begin{equation}\label{orasoon}
\Gamma^{k(s)}_{ij}(1)=0\ \ \textup{for all}\ \ i,j,k\in \{1,\ldots,n\}.\end{equation}
For every 
 non-degenerate geodesic $c\in R^n$ we have
 $$\delta^{(s)}v_k=0 \ \ \textup{for all}\ \ k\in \{1,\ldots,n\};$$
 equivalently, each $v_k$ is either $0$ or a root of unity.
 This should be thought of as saying that the geodesics, in this case, are `straight lines.'
 For conditions implying (\ref{orasoon}) see Remark \ref{scholl} below.
 Similarly if $\lambda^{(0)}:=c(P_0)$ and
 $$(\Gamma^{k(s)}_{ij}(1))(\lambda^{(0)})=0\ \ \textup{for all}\ \ i,j,k\in \{1,\ldots,n\}$$
 then, in view of formula (\ref{lowonbread}),
 $$\delta^{(s)}v_k\equiv 0 \ \ \textup{mod $p$ for all}\ \ k\in \{1,\ldots,n\}$$
 i.e., the `acceleration' of $c$ vanishes at $P_0$.
 For a partial converse of this see Remark \ref{piinee}.\end{remark}

By Propositions \ref{eqforgeo} and \ref{flowww} we get the following results:

\begin{theorem}\label{eupar} (Parallel transport). 
Let $\Delta^{(s)G}$ be a $\pi$ -connection of degree $s$ on $G$.  For every non-degenerate curve 
 $c\in R^n$ and every vector $w^{(0)}\in R^n$
there exists a unique vector  $w\in R^n$ such that $w$ is parallel along $c$ 
with respect to
$\Delta^{(s)G}$
and $$w\equiv w^{(0)}\ \ \textup{mod}\ \ p.$$
\end{theorem}

\begin{theorem}\label{eug} (Existence and uniqueness of geodesics). 
Let $\Delta^{(s)G}$ be a $\pi$ -connection of degree $s$ on $G$.  For every non-degenerate curve 
 $c^{(0)}\in R^n$ 
there exists a unique non-degenerate curve  $c\in R^n$ such that $c$ is a geodesic for 
$\Delta^{(s)G}$
and $$\begin{array}{rcll}
c &\equiv & c^{(0)} & \textup{mod}\ \ p,\\
\ & \ & \ & \ \\
\delta^{(s)}c & \equiv & \delta^{(s)}c^{(0)} & \textup{mod}\ \ p.\end{array}$$
\end{theorem}

\begin{remark}
Let $\Delta^{(s)G}$ be a $\pi$-connection of degree $s$ on $G$.
Using the concept of evaluation maps $R\rightarrow k$, $a\mapsto a(P)$, cf. Definition \ref{evaeva}, one can interpret the Theorem~\ref{eupar} and Theorem~\ref{eug} as follows. Recall  $P_0=t\in R\{t\}$ and fix  $P\in R\{t\}$.

1) For every non-degenerate curve  $c\in R^n$  we can define the {\bf parallel transport map}  along $c$ from $P_0$ to $P$,
$$\textup{par}_{c,P}:k^n\rightarrow k^n,$$
 as follows. Let $\lambda\in k^n$, consider the unique vector  $w\in R^n$ parallel along $c$ 
and satisfying 
$$\lambda=w(P_0);$$ then set 
$$\textup{par}_{c,P}(\lambda):=w(P).$$
Note that, by Proposition \ref{regularr} the map $\textup{par}_{c,P}$ is induced by a morphism of schemes over $k$, which we still denote by
$$\textup{par}_{c,P}:\mathbb A_k^n\rightarrow \mathbb A_k^n.$$
Furthermore let
 $$P=\sum_{d=D_0}^D P_{(d)}\in R\{t\}$$
 with $P_{(d)}$ homogeneous of degree $d$ and $D_0\leq D$.  For every solution $w$ of the system 
 (\ref{pararara}) and every root of unity $\zeta\in R$ we recall from   Remark \ref{syyster}, assertion (1)
 that $\zeta w$ is still a solution of the system. On the other hand 
 we have
 $$\begin{array}{rcl}
 (\zeta  w_i)(P_0) &= & \zeta  w_i(P_0),\\
 \ & \ & \ \\
 (\zeta w_i)(P) & = & P^{(s)}\{\zeta w_i\} (P_0)\\
 \ & \ & \ \\
 \ & = & \sum_{d=D_0}^D P_{(d)}^{(s)}\{\zeta w_i\}(P_0)\\
 \ & \ & \ \\
 \ & = & \sum_{d=D_0}^D \zeta^d  P_{(d)}^{(s)}\{w_i\}(P_0).\end{array}$$
 So we have  the following formula:
 $$\textup{par}_{c,P}(\gamma\cdot \lambda)=\sum_{d=D_0}^D \gamma^d\cdot \textup{par}_{c,P_{(d)}}(\lambda),\ \ \textup{for all}\ \ \gamma\in k, \ \lambda\in k^n.$$
 In particular $\textup{par}_{c,P}$ sends the lines of $k^n$ passing through the origin into
 curves in $k^n$ parameterized by polynomials of degree $\leq D-D_0+1$.
 
2)  Recalling that $p/\pi=\sum r_i\theta_i\in R_{\pi}$, $r_i\in R$,
let $H_0\subset k^n$ be the hyperplane of all vectors orthogonal  to the vector $(r_1,\ldots,r_n)$ 
and let $H_1\subset k^n$ be the hyperplane of all vectors orthogonal to $(1,\ldots,1)$.
We can define the {\bf geodesic exponential map} 
$$\textup{exp}:(k^n\setminus H_0)\times (k^n\setminus H_1)\rightarrow k^n,\ \ \ (\lambda_0,\lambda_1)\mapsto \textup{exp}_{\lambda_0}(\lambda_1),$$
as follows. Let $\lambda_i\in k^n\setminus H_i$ and consider the unique geodesic  $c\in R^n$ such that
 $$\lambda_0=c(P_0)\ \ \textup{and}\ \ \lambda_1=v(P_0)$$
 where $v$ is the velocity of $c$;
  then set 
$$\textup{exp}_{\lambda_0}(\lambda_1):=c(P).$$
Note that, intuitively,  the `origins' of our geodesics are assumed to avoid the hyperplane $H_0$ while the directions of our geodesics at their common origin are assumed to avoid the hyperplane $H_1$. 
Again, by Proposition \ref{regularr}, the map of sets $\textup{exp}$ is induced by a morphism of schemes, still denoted by
$$\textup{exp}:(\mathbb A_k^n\setminus \mathbb H_0)\times (\mathbb A_k^n\setminus \mathbb H_1)\rightarrow \mathbb A_k^n,$$
where $\mathbb H_i\subset \mathbb A^n$ are the closed subschemes corresponding to our hyperplanes $H_i$.
\end{remark}

 \subsection{Torsion symbols and symmetric $\pi$-connections}

\begin{definition} \label{deftorsym}
Let $y=(y_1,\ldots,y_n)$ be an $n$-tuple of matrices of $n\times n$ indeterminates
and consider the polynomial ring
$$R_{\pi}[y]:=R_{\pi}[(y_i)_{jk}\ |\ i,j,k\in \{1,\ldots,n\}].$$  
A {\bf torsion symbol of the second kind} is 
an $n$-tuple $L=(L^1,\ldots,L^n)$  of antisymmetric matrices $L^k\in 
\textup{Mat}_n(\widehat{R_{\pi}[y]})$,
$$L_{ij}^k+L_{ji}^k=0.$$ 
\end{definition}

We will  denote by $L^{(s)}=(L^{1(s)},\ldots,L^{n(s)})$ the torsion symbols that
appear in the context of higher $\pi$-derivations of degree $s$. 

Recall the notation $1=1_n$ for the identity matrix. From Notation~\ref{notynoty}, we denoted the natural substitution map $\mathcal A \to \mathcal A/\mathcal P \cong R_\pi$ where $\mathcal A = R_\pi[x,\det(x)^{-1}]$ by the notation $F(1)$ for $F \in \mathcal A$. Similarly, given matrices $\Xi = (\Xi_1,\ldots,\Xi_n)$ in $\textup{Mat}_n(S)$ for an $R_\pi$-algebra $S$, we denote the natural substitution $\widehat{R_{\pi}[y]} \to S$, replacing $y_i$ by $\Xi_i$ as $F(\Xi) := F(\Xi_1,\ldots,\Xi_n)$ for $F \in \widehat{R_\pi[y]}$. When $\Xi_i = 1$ for all $i$, this is denoted again simply by $F(1)$. Specifically, this process defines for the matrices $\Lambda^{(s)} = (\Lambda_1^{(s)},\ldots,\Lambda_n^{(s)})$, $\Lambda_i^{(s)}:=(x^{(p^s)})^{-1}(\phi_i^{(s)})^G (x) \in \textup{GL}_N(\mathcal A)$, elements $L_{ij}^{k(s)}(\Lambda^{(s)}) \in \mathcal A$ for a torsion symbol $L^{(s)}$. 


\begin{remark}
For now we do not impose any restriction on the  torsion symbol. In Part 2 of this paper we will show how to canonically choose torsion symbols  that `measure the non-commutation' of $\phi^{(s)}_1,\ldots,\phi_n^{(s)}$ in $R_{\pi}$.
\end{remark}

\begin{remark}\label{zoiu}
The following trivial fact will play a role later. Let $S$ be a ring, let $X=(X_{ij})\in \textup{GL}_n(S)$ be  an invertible matrix, let $Y_i=((Y_i)_{jk})\in \textup{Mat}_n(S)$ be matrices, and let 
$Z^k=(Z_{ij}^k)\in \textup{Mat}_n(S)$ be matrices, where $i,j,k\in \{1,\ldots, n\}$. The  following are equivalent:

\smallskip

1) $(Y_i)_{kj}-(Y_j)_{ki}=Z^k_{ij}$ for all $i,j,k$.

\smallskip

2) $(XY_i)_{kj}-(XY_j)_{ki}=\sum_m Z^m_{ij}X_{km}$ for all $i,j,k$.

\end{remark}

\begin{definition} \label{defsym} Assume $N=n$.
 A  $\pi$-connection  of degree $s$ as in Definition \ref{defconnn}  is called {\bf symmetric} with respect  to a torsion symbol of the second kind
$L^{(s)}=(L^{1(s)},\ldots,L^{n(s)})$ if for all $i,j,k\in \{1,\ldots,n\}$ we have
\begin{equation}\label{uff22}
\Gamma^{k(s)}_{ij}-\Gamma^{k(s)}_{ji}= L^{k(s)}_{ij}(\Lambda^{(s)}).\end{equation}
\end{definition}

\begin{remark}\label{piinee}
Assume  $L^{(s)}(1)\equiv 0$ mod $\pi$. Natural classes of examples  when this happens will be considered in Part 2.
Assume furthermore that we are given a $\pi$-connection that is symmetric with respect to $L^{(s)}$ and assume we have a point $\lambda^{(0)}=(\lambda^{(0)}_1,\ldots,\lambda^{(0)}_n)\in k^n$ such that $\sum_{i=1}^nr_i\lambda^{(0)}_i\neq 0$  where we recall that $\pi/p = \sum r_i \theta_i \in R_\pi$. Assume that all non-degenerate geodesics $c$ with $c(P_0)=\lambda^{(0)}$ and velocity $v$ 
 satisfy $\delta v\equiv 0$ mod $p$.
Then the Christoffel symbols vanish at $\lambda^{(0)}$:
$$(\Gamma^{k(s)}_{ij}(1))(\lambda^{(0)})=0.$$
This is analogous to the well known similar statement in classical differential geometry in the following sense.
By formula (\ref{geeoo})  we get that
$$\sum_{ij}(c^*(\Gamma^{k(s)}_{ij}))(P_0)\xi_i\xi_j=0$$
for all $k$ and all vectors $\xi:=(\xi_1,\ldots,\xi_n)\in k^n$ with $\sum_{i=1}^n\xi_i\neq 0$. We deduce that
the same holds for all vectors $\xi\in k^n$. By formula  (\ref{lowonbread}) we get
$$\sum_{ij}\Gamma^{k(s)}_{ij}(\lambda^{(0)})\xi_i\xi_j=0$$
for all $\xi\in k^n$. Hence
$$\Gamma^{k(s)}_{ij}(\lambda^{(0)})+\Gamma^{k(s)}_{ji}(\lambda^{(0)})=0$$
and we conclude by symmetry.
\end{remark}

\subsection{Metric $\pi$-connections}

\begin{definition}
 By a {\bf metric} of dimension $N$, over $R_{\pi}$ we mean a symmetric matrix $q\in \textup{GL}_N(R_{\pi})^{\textup{sym}}$.\end{definition}
 
 We sometimes denote by $q^{(s)}$ the metrics that
appear in the context of higher $\pi$-derivations of degree $s$.
 We view $q$ as an arithmetic analogue of a metric on a vector bundle of rank $N$ 
 over a manifold.
 When $N=n$ we may also view $q$ as an arithmetic analogue of a metric on the tangent bundle.

\begin{definition}\label{qui1}
For a metric $q^{(s)}$ we define the $R_{\pi}$-algebra map,
$$\cH_{q^{(s)}}:\mathcal A\rightarrow \mathcal A,\ \ \cH_{q^{(s)}(x)}:=x^tq^{(s)}x,$$
i.e.,
$\cH_{q^{(s)}}(x_{ij})=\sum_{k,l}q^{(s)}_{kl}x_{ki}x_{lj}$.    
 A  $\pi$-connection  of degree $s$  as in Definition \ref{defconnn} is said to be {\bf metric} with respect to  $q^{(s)}$ 
  if the following diagrams are commutative for $i\in \{1,\ldots,n\}$:
 \begin{equation}
 \label{got}
 \begin{array}{rcl}
\mathcal A& \stackrel{(\phi_{i,0}^{(s)})^G}{\longrightarrow} & \mathcal A \\
 \cH_{q^{(s)}} \downarrow &\ &\downarrow \cH_{q^{(s)}}\\
 \mathcal A & \stackrel{(\phi_i^{(s)})^G}{\longrightarrow} & \mathcal A\end{array}
\end{equation}
\end{definition}

Here we recall that $((\phi_{i,0}^{(s)})^G)$ corresponds to the trivial  $\pi$-connection of degree $s$; see Definition \ref{tritri}.
For a discussion of the analogy with classical differential geometry see the Appendix of Part 2 of this paper.

\subsection{Formulae}\label{sst}

 Explicitly for every 
 $q^{(s)}\in \textup{GL}_N(R_{\pi})^{\textup{sym}}$ and every  $\pi$-connection $\Delta^{(s)G}=((\delta_{1}^{(s)})^G,\ldots,(\delta_{n}^{(s)})^G)$  of degree $s$ 
we consider the following $N\times N$ matrices 
 with entries in $\mathcal A$:
\begin{equation}
\label{1955}
\begin{array}{rcl}
A_i^{(s)} & := & x^{(p^s)t}\cdot \phi^{(s)}_{i}(q^{(s)})\cdot x^{(p^s)},\\
\ & \ & \ \\
 B^{(s)} & := & (x^tq^{(s)}x)^{(p^s)}\end{array}\end{equation}

\begin{lemma}\label{eajuns}
If $C_i^{(s)}:=\frac{1}{\pi}(B^{(s)}-A_i^{(s)})$ then
\begin{equation}
C^{(s)}_i\equiv -\delta_i^{(s)}q^{(s)}\ \ \textup{mod}\ \ \mathcal P\end{equation}
in the ring $\mathcal A$.\end{lemma}

Here we recall that $\mathcal P:=(x-1)$ is the ideal generated by the entries $x_{kl}-\delta_{kl}$ of the matrix $x-1$.

\medskip

{\it Proof}.  We have that 
$\pi C_i^{(s)}+\pi\delta_i^{(s)} q^{(s)}$ has entries in the ideal
$(\pi)\cap \mathcal P$. But the latter intersection ideal equals the product ideal $(\pi)\mathcal P$ because $\mathcal P$ is a prime ideal not containing $\pi$. So $\pi C^{(s)}_i+\pi \delta_i^{(s)} q^{(s)}=\pi M_i$ with $M_i$ a matrix with entries in $\mathcal P$ hence $C^{(s)}_i+\delta_i^{(s)} q^{(s)}=M_i$. \qed

\begin{remark}\label{comua} \ 

1) For $\mathcal M = (\pi,\mathcal P)$, the reduction modulo $\mathcal M^2$ of $\phi_i^{(s)}$ is uniquely determined by 
the reduction modulo  $\mathcal M^2$ of $\Lambda_i^{(s)}$; the latter is completely determined by the reduction modulo $\mathcal M$ of $\Gamma_i^{(s)}$.

2) The commutativity of \eqref{got} is  equivalent to the matrix equalities:
\begin{equation}
\label{LALB}
\Lambda_{i}^{(s)t}A_i^{(s)} \Lambda_{i}^{(s)}=B^{(s)}.\end{equation}

3) The symmetry condition (\ref{uff22}) is equivalent to the equalities:
\begin{equation}
\label{uf22}
(\Lambda_i^{(s)}-1)_{kj}-(\Lambda_j^{(s)}-1)_{ki}=\pi\cdot L_{ij}^{k(s)}(\Lambda^{(s)}).
\end{equation}
\end{remark}

\begin{notation}
We consider  the  `lowering of the indices' operation by defining
\begin{equation}
\label{loweringg}
\begin{array}{rcl}
L_{ijk}^{(s)} & := & \sum_m L_{ij}^{m(s)}(q_{mk}^{(s)})^{p^s},\\
\ & \ & \ \\
\Gamma_{ijk}^{(s)} & := & \sum_m \Gamma_{ij}^{m(s)} (q_{mk}^{(s)})^{p^s}.\end{array}\end{equation}
Note that, consistent with our previously explained notation,
 $\Gamma_{ijk}^{(s)}$ is {\it not} the $jk$ entry of $\Gamma_i^{(s)}$, the latter being denoted by $(\Gamma_i^{(s)})_{jk}=\Gamma^{k(s)}_{ij}$, see Equation \ref{trematt}. A similar remark holds 
 for $L_{ijk}^{(s)}$.\end{notation}

\begin{definition}\label{firstkindd}
$\Gamma_{ijk}^{(s)}$ are the {\bf Christoffel symbols of the first kind} (relative to the metric $q^{(s)}$);
$L_{ijk}^{(s)}$ are the {\bf torsion symbols of the first kind} (relative to the metric $q^{(s)}$).
\end{definition}

\subsection{Arithmetic Levi-Civita connection}

Recall that $\mathcal P$ is the prime ideal in $\mathcal A$
generated by  the entries of the matrix $x-1$ and we consider the maximal ideal
 $\mathcal M:=(\pi,\mathcal P)$ of $\mathcal A$.
Our main result in this subsection is the following arithmetic PDE analogue of the Fundamental Theorem of Riemannian Geometry. 

\begin{theorem}\label{LCC}
(Arithmetic Levi-Civita Theorem)
Assume $N=n$.
Assume $s\geq 1$ is an integer, $q^{(s)}\in G(R_{\pi})^{\textup{sym}}$ is a metric,   and $L^{(s)}=(L^{1(s)},\ldots,L^{n(s)})$ is a torsion symbol of the second kind. 

\medskip

1) There exists a unique   $\pi$-connection $\Delta^{(s)\textup{LC}}=(\delta_1^{(s)\textup{LC}}, \ldots, \delta_n^{(s)\textup{LC}})$  of degree $s$  on $G$ that is metric with respect to $q^{(s)}$ and symmetric with respect to $L^{(s)}$.

\medskip

2) If $(\Gamma_{ijk}^{\textup{(s)\textup{LC}}})$ and $(L_{ijk}^{(s)})$ are  the corresponding Christoffel symbols and the torsion symbols of the first kind
 then we have the following congruence in $\mathcal A$:
 
 \medskip

$$\Gamma_{ijk}^{(s)\textup{LC}} \equiv -\frac{1}{2}(\delta_{i}^{(s)}q^{(s)}_{jk}+\delta_{j}^{(s)}q^{(s)}_{ik}-\delta_{k}^{(s)}q^{(s)}_{ij})+
\frac{1}{2}
(L^{(s)}_{kij}(1)+L^{(s)}_{ijk}(1)-L^{(s)}_{jki}(1))
\ \ \ \textup{mod} \ \ \ \mathcal M.$$

\medskip

3) If $m\geq 1$ is an integer
such  that for all $i,j,k\in \{1,\ldots,n\}$ we have the congruences

\medskip

$$L_{ij}^{k(s)}(1)\equiv \delta_i^{(s)}q_{jk}^{(s)}\equiv 0\ \ \ \textup{mod}\ \ \pi^{m}$$

\medskip

\noindent in $R_{\pi}$ then  we have the following congruences in $\mathcal A$:

 \medskip
 
 $$\Gamma_{ijk}^{(s)\textup{LC}}\equiv 0 \ \ \textup{mod}\ \ (\pi^{m},\mathcal P).$$

\end{theorem}

 \begin{remark}
 Assume that  $L^{(s)}_{ijk}(1)\equiv 0$ mod $\pi$ for all $i,j,k$.
 The  equation in Theorem \ref{LCC} (2)  reads
  \begin{equation}
  \label{mict}
  \Gamma_{ijk}^{(s)\textup{LC}} \equiv -\frac{1}{2}(\delta_{i}^{(s)}q^{(s)}_{jk}+\delta_{j}^{(s)}q^{(s)}_{ik}-\delta_{k}^{(s)}q^{(s)}_{ij})
\ \ \ \textup{mod} \ \ \ \mathcal M.
  \end{equation}
  The latter equation is of course analogous to the classical expression of the Christoffel symbols with respect to the coordinate vector fields.\end{remark}
    
  \begin{remark}\label{scholl}
 \ 
 
 1) If $L_{ij}^{k(s)}(1)=\delta_i^{(s)}q_{jk}^{(s)}=0$ for all $i,j,k$ then $\Gamma^{k(s)}_{ij}(1)=\Gamma_{ijk}^{(s)\textup{LC}}(1)=0$. 
  
  2) The condition $\delta_i^{(s)}q_{jk}^{(s)}=0$ is satisfied for some $i,j,k$ if and only if 
  $q_{jk}^{(s)}$ is either $0$ or a root of unity in $R$; cf. Lemma \ref{vinne}.
  
  3) Natural examples when $L_{ij}^{k(s)}(1)=0$ for all $i,j,k$ will be discussed in Part 2 of this paper.
  \end{remark}

\begin{remark}
We take the opportunity here to correct a typo in \cite{Bu17}: in each of the Equations 4.88 and 4.89 on page 197 in loc.cit. the signs $+$ and $-$ need to be interchanged.
\end{remark}

\begin{definition}\label{arlecico}
The unique  $\pi$-connection  $\Delta^{(s)\textup{LC}}$  of degree $s$ in Theorem \ref{LCC} is called the {\bf arithmetic  Levi-Civita connection} of degree $s$ attached to $q^{(s)}$ and $L^{(s)}$. \end{definition}

Our Theorem \ref{LCC} will follow  from the following construction, which is an enhancement of \cite[Thm. 2.15]{Bu19}.

\begin{theorem}\label{theone}
(Abstract Arithmetic Levi-Civita Theorem)
Let $S$ be a ring in which $2$ is invertible, $P$ an ideal in $S$, $\pi\in S$ a regular element in $S$ whose image in $S/P$ is also a regular element, and assume $S$ is $\pi$-adically complete. Let $n\geq 1$ and $m\geq 0$ be  integers. Consider matrices $A_i,B\in \textup{GL}_n(S)$, $i\in \{1,\ldots,n\}$. Moreover consider an $n$-tuple $y=(y_1,\ldots,y_n)$  of $n\times n$ matrix indeterminates $y_1,\ldots,y_n$ and a tuple $L:=(L^1,\ldots,L^n)$, $L^k\in \textup{Mat}_n(\widehat{S[y]})$. Assume the following properties:

\smallskip

1) $A_i$ and $B$ are symmetric, i.e., $(A_i)_{jk}=(A_i)_{kj}$, $B_{jk}=B_{kj}$;

\smallskip

2)  $L^k$ are antisymmetric, i.e., $L_{ij}^k+L_{ij}^k=0$;

\smallskip

3)  $A_i\equiv B$ mod $\pi$;

\smallskip

4) $L^k\equiv 0$ mod $(\pi^m,P)$ and $A_i\equiv B$ mod $(\pi^{m+1},P)$.

\smallskip

\noindent Then there exists a unique $n$-tuple $\Lambda=(\Lambda_1,\ldots,\Lambda_n)$ of matrices $\Lambda_i\in \textup{GL}_n(S)$
 satisfying the following properties:

\smallskip

a) $\Lambda_i\equiv 1$ mod $\pi$;

\smallskip

b) $\Lambda_i^t A_i \Lambda_i=B$;

\smallskip

c)  $(\Lambda_i-1)_{kj}-(\Lambda_j-1)_{ki}=\pi L_{ij}^k(\Lambda)$;

\smallskip

d) $\Lambda_i\equiv 1$ mod $(\pi^{m+1},P)$.

\smallskip

\noindent Moreover if 
$$C_i:=\frac{1}{\pi}(B-A_i),\ \ L_{ijk}:=\sum_l L_{ij}^lB_{kl}$$
 then the following congruences hold in $S$:
 
 \smallskip
 
$$
\frac{1}{\pi}(B(\Lambda_i-1))_{kj}\equiv \frac{1}{2}((C_i)_{jk}+(C_j)_{ik}-(C_k)_{ij})+\frac{1}{2}
(L_{kij}(1)+L_{ijk}(1)-L_{jki}(1))\ \ \textup{mod}\ \ \pi.$$
\end{theorem}

\medskip

For $m=0$ condition (4) above is a consequence of condition (3); and similarly condition $(d)$ is a consequence of condition $(a)$. So for $m=0$ conditions (4) and $(d)$ can be deleted from the statement of the Theorem.

\medskip

{\it Proof}. The proof runs along the lines of \cite[Thm. 2.15]{Bu19}  but it is  
more 
involved 
and it contains  new tricks (related to the matrices $F_i$ below). 

 We will   construct by induction a sequence of $n$-tuples
\begin{equation}
\label{lamb}
\Lambda^{\nu}=(\Lambda_1^{\nu},\ldots,\Lambda_n^{\nu}),\ \ \nu\geq 1\end{equation}
of $n\times n$ matrices with entries in $S$
such that for $\nu\geq 1$ and all $i,j\in \{1,\ldots,n\}$
 the following properties hold:

\medskip

i) $\Lambda_i^{1}=1$;

\smallskip

ii) $\Lambda_i^{\nu+1}\equiv \Lambda_i^{\nu}$ mod $\pi^{\nu}$;

\smallskip

iii) $\Lambda_i^{\nu t}A_i\Lambda_i^{\nu}\equiv B$ mod $\pi^{\nu}$;

\smallskip

iv) $(B(\Lambda_i^{\nu}-1))_{kj}-(B(\Lambda_j^{\nu}-1))_{ki}\equiv \pi L_{ijk}(\Lambda^{\nu})\ \ \textup{mod}\ \ \pi^{\nu}$;

\smallskip

v) $\Lambda_i^{\nu}\equiv 1$ mod $(\pi^{m+1},P)$ for $\nu\leq m+1$.

\smallskip

\noindent  Here the superscript $\nu$ of $\Lambda_i$ (and of other matrices to be considered below) is an index (not an exponent).
Taking $\Lambda_i$ to be the limit of $\Lambda_i^{\nu}$ as $\nu\rightarrow \infty$  and using Remark \ref{zoiu} ends the proof.

To construct  our sequence of $n$-tuples (\ref{lamb}) define the $n$-tuple for $\nu=1$ by condition (i) and note that  conditions (ii,) (iii), (iv) are then automatically satisfied.
 Assume now the $n$-tuple (\ref{lamb})
was constructed for some $\nu\geq 1$ and seek the $n$-tuple (\ref{lamb})
corresponding to $\nu+1$ in the form
\begin{equation}
\label{alpha}
\Lambda_i^{\nu+1}=\Lambda_i^{\nu}+\pi^{\nu}Z^{\nu}_i\end{equation}
where $Z^{\nu}_i$ needs to be determined.
Write 
\begin{equation}
\label{betabeta1}
\Lambda_i^{\nu t}A_i\Lambda_i^{\nu}= B-\pi^{\nu}C^{\nu}_i,\end{equation}
\begin{equation}
\label{betabeta2}
(B(\Lambda_i^{\nu}-1))_{kj}-(B(\Lambda_j^{\nu}-1))_{ki}= \pi L_{ijk}(\Lambda^{\nu})-\pi^{\nu}L_{ijk}^{\nu},
\end{equation}
for matrices $C^{\nu}_i$ and $L_{ijk}^{\nu}$ with entries in $S$.
In particular,
\begin{equation}
\label{inca1}
L_{ijk}^{1}=L_{ijk}(1).\end{equation}
Also, using the fact that $\pi$ is  a regular element in $S$ and $S/P$,  we have that
\begin{equation}\label{aldoo}
C^{\nu}_i\equiv L_{ijk}^{\nu}\equiv 0\ \ \textup{mod}\ \ \pi^{m+1-\nu}\ \ \textup{if}\ \ \nu\leq m+1.
\end{equation}
By (\ref{alpha}) and (\ref{betabeta1}) we get
\begin{equation}
\label{cip cirip}
\begin{array}{rcl}
\Lambda_i^{(\nu+1)t}A_i\Lambda_i^{\nu+1} & \equiv & \Lambda_i^{\nu t}A_i\Lambda_i^{\nu} +\pi^{\nu}(\Lambda_i^{\nu t}A_iZ^{\nu}_i+Z^{\nu t}_iA_i\Lambda_i^{\nu})\ \ \text{mod}\ \ \pi^{\nu+1}\\
\  & \  & \  \\
\  & \equiv & B+\pi^{\nu}(-C^{\nu}_i+BZ^{\nu}_i+Z^{\nu t}_iB) \ \ \text{mod}\ \ \pi^{\nu+1}.
\end{array}\end{equation}
Now 
$A_i^t=A_i$ and $B^t=B$ so 
$C_i^{\nu t}=C^{\nu}_i.$
 Hence $(C^{\nu}_i)_{jk}=(C^{\nu}_i)_{kj}$.
 Define the matrices
 $D^{\nu}_i=((D^{\nu}_i)_{jk})$
by setting
\begin{equation}
\label{gamma}
(D^{\nu}_i)_{jk}:=\frac{1}{2}((C^{\nu}_i)_{jk}+(C^{\nu}_j)_{ik}-(C^{\nu}_k)_{ij}).\end{equation} 
Then
\begin{equation}
\label{66}
(D^{\nu}_i)_{jk}=(D^{\nu}_j)_{ik}\end{equation}
 and $(D^{\nu}_i)_{jk}+(D^{\nu}_i)_{kj}=(C^{\nu}_i)_{jk}.$
So we have:
\begin{equation}
\label{tinna}
D^{\nu}_i+D_i^{\nu t}=C^{\nu}_i.\end{equation}
On the other hand, by (\ref{betabeta2}), we get $L^{\nu}_{ijk}+L^{\nu}_{jik}=0.$

Define the matrices
$$F^{\nu}_i=((F^{\nu}_i)_{jk})$$
by setting
\begin{equation*}
(F^{\nu}_i)_{jk}:=\frac{1}{2}(L^{\nu}_{kij}+L^{\nu}_{ijk}-L^{\nu}_{jki}).
\end{equation*}
So by (\ref{inca1}),
\begin{equation}
\label{monster1}
(F^{1}_i)_{jk}= \frac{1}{2}(L_{kij}(1)+L_{ijk}(1)-L_{jki}(1)).
\end{equation}
For all $\nu\geq 1$ one has
\begin{equation}
\label{777}
(F^{\nu}_i)_{jk}-(F^{\nu}_j)_{ik}=L^{\nu}_{ijk},
\end{equation}
\begin{equation*}
(F^{\nu}_i)_{jk}+(F^{\nu}_i)_{kj}=0;
\end{equation*}
hence
\begin{equation}
\label{feyy}
F^{\nu}_i+F^{\nu t}_i=0.\end{equation}
Setting 
\begin{equation}
\label{omega}
Z^{\nu}_i:=B^{-1}(D^{\nu t}_i+F_i^{\nu t})\end{equation}
  we have
 \begin{equation}
 \label{dorra1}
 BZ^{\nu}_i=D^{\nu t}_i+F_i^{\nu t},\end{equation}
 \begin{equation*}
 Z^{\nu t}_iB=D^{\nu}_i+F_i^{\nu}.\end{equation*}
 Moreover using (\ref{aldoo}) we have
 $Z^{\nu}_i\equiv 0\ \ \textup{mod}\ \ (\pi^{m+1}, P)\ \ \textup{for}\ \ \nu\leq m+1,$
 hence, by (\ref{alpha}), we have
 $  
 \Lambda^{\nu+1}_i\equiv 1\ \ \textup{mod}\ \ (\pi^{m+1},P)\ \ \textup{for}\ \ \nu\leq m+1,
 $
 which shows the condition (v) holds for $\Lambda_i^{\nu+1}$.
 On the other hand,  using (\ref{tinna}) and (\ref{feyy}), we get
 \begin{equation}
 \label{dorra3}
 BZ^{\nu}_i+Z^{\nu t}_iB=C_i^{\nu}.\end{equation}
 By (\ref{cip cirip}) and (\ref{dorra3}) we get  
 $$\Lambda_i^{(\nu+1)t}A_i\Lambda_i^{\nu+1}\equiv B\ \ \ \text{mod}\ \ \pi^{\nu+1}$$
 and hence condition (iii) holds for $\Lambda_i^{\nu+1}$ mod $\pi^{\nu+1}$.

 To check condition (iv) for $\nu$ replaced by $\nu+1$ note that 
 for 
 $$\lambda_{ikj}:=(B(\Lambda_i^{\nu+1}-1))_{kj}$$
 we have
 $$
 \begin{array}{rcl}
\lambda_{ikj} & = & (B(\Lambda_i^{\nu}+\pi^{\nu}Z^{\nu}_i-1))_{kj}\ \ \textup{cf}. (\ref{alpha})\\
\  & \  & \  \\
 \ &  = & (B(\Lambda^{\nu}_i-1))_{kj}+\pi^{\nu}(D_i^{\nu})_{jk}+\pi^{\nu} (F^{\nu}_i)_{jk}\ \ \textup{cf}.\ (\ref{dorra1})   \\
\  & \  & \  \\
\  & = & (B(\Lambda^{\nu}_i-1))_{kj}+\pi^{\nu}(D_j^{\nu})_{ik}+\pi^{\nu} (F^{\nu}_j)_{ik}+\pi^{\nu}L_{ijk}^{\nu}\ \ \textup{cf}. \ (\ref{66}),(\ref{777})\\
\  & \  & \  \\
\  & = & (B(\Lambda^{\nu}_j-1))_{ki}+\pi^{\nu}(D_j^{\nu})_{ik}+\pi^{\nu} (F^{\nu}_j)_{ik}+\pi L_{ijk}
(\Lambda^{\nu})\ \ \textup{cf}. \ (\ref{betabeta2})\\
\  & \  & \  \\
\  & \equiv & (B(\Lambda^{\nu}_j-1))_{ki}+\pi^{\nu}(D_j^{\nu})_{ik}+\pi^{\nu} (F^{\nu}_j)_{ik}+
\pi L_{ijk}
(\Lambda^{\nu+1})\ \ \textup{mod}\ \ \pi^{\nu+1}\\
\  & \  & \  \\
\  & = & (B(\Lambda_j^{\nu}+\pi^{\nu}Z^{\nu}_j-1))_{ki}+\pi L_{ijk}(\Lambda^{\nu+1})\ \ \textup{cf}. \ (\ref{dorra1})\\
\  & \  & \  \\
\  & = & \lambda_{jki}+\pi L_{ijk}(\Lambda^{\nu+1})\ \ \textup{cf}. \ (\ref{alpha}).
\end{array}
$$
So condition (iv) holds for $\nu$ replaced by $\nu+1$. This ends the proof of the existence part of our theorem.

 We next prove the uniqueness part 
  in our theorem.
  Assume we have two tuples which we denote by 
  $$\Lambda=(\Lambda_1,...,\Lambda_n)\ \ \ \text{and}\ \ \  \Lambda'=(\Lambda'_1,...,\Lambda'_n)$$
 such that for al $i,j,k$ we have
 \begin{equation*}
 \Lambda_i\equiv \Lambda_i'\equiv 1\ \ \textup{mod}\ \ \pi,
 \end{equation*}
  \begin{equation}
  \label{fish}
  \Lambda_i^tA_i\Lambda_i= (\Lambda_i')^tA_i\Lambda'_i=B,\end{equation}
  \begin{equation}
  \label{clouds}
  (B(\Lambda_i-1))_{kj}-(B(\Lambda_j-1))_{ki}=\pi L_{ijk}(\Lambda),\end{equation}
  \begin{equation}
  \label{clouds'}
  (B(\Lambda'_i-1))_{kj}-(B(\Lambda'_j-1))_{ki}=\pi L_{ijk}(\Lambda').
  \end{equation}
  We will prove that for all $\nu\geq 1$ we have
  \begin{equation}
  \label{ciri}
  \Lambda_i\equiv \Lambda'_i\ \ \ \text{mod}\ \ \pi^{\nu}
  \end{equation}
  and this will end the proof of the uniqueness part. We proceed 
   by induction on $\nu$. The case $\nu=1$ is clear. Assume \ref{ciri} holds for some $\nu\geq 1$ and write
   \begin{equation}
   \label{toctoc}
   \Lambda'_i=\Lambda_i+\pi^{\nu}Z^{\nu}_i.\end{equation}
   From (\ref{fish}) we get
   $$B\equiv B+\pi^{\nu}A_iZ^{\nu}_i+\pi^{\nu}Z^{\nu t}_iA_i\equiv
   B+\pi^{\nu}BZ^{\nu}_i+\pi^{\nu}Z^{\nu t}_iB
   \ \ \ \text{mod}\ \ \ \pi^{\nu+1},
   $$
   hence, setting
   $$E_i:=Z^{\nu t}_iB=:((E_i)_{jk})$$
   we get
   $E_i+E_i^t\equiv 0\ \text{mod}\ \pi,$
   hence
   \begin{equation}
   \label{tuctuc}
   (E_i)_{jk}\equiv - (E_i)_{kj}\ \ \ \text{mod}\ \ \ \pi.
   \end{equation}
   On the other hand,  subtracting (\ref{clouds}) from (\ref{clouds'}) and using (\ref{toctoc})
   we get
   \begin{equation}
   \label{tctc}
   (E_i)_{jk}\equiv (E_j)_{ik}\ \ \textup{mod}\ \ \pi.
   \end{equation}
   Combining \ref{tuctuc} and \ref{tctc} we get
   \begin{equation}
   \label{tirgo}
   (E_i)_{jk}\equiv (E_j)_{ik}\equiv - (E_j)_{ki}\ \ \ \text{mod}\ \ \ \pi.\end{equation}
   Applying \ref{tirgo} three times we get
   $$(E_i)_{jk}\equiv - (E_j)_{ki}\equiv (E_k)_{ij}\equiv - (E_i)_{jk}\ \ \ \text{mod}\ \ \pi,$$
 hence $2(E_i)_{jk}\equiv 0\ \text{mod}\ \pi.$
 Since $2$ is invertible in $S$ we get $(E_i)_{jk}\equiv 0\ \text{mod}\ \pi$
hence
  $Z^{\nu}_i\equiv 0\ \text{mod}\ \pi,$ which implies
   $\Lambda'_i\equiv \Lambda_i\ \text{mod}\ \pi^{\nu+1},$
   and our induction step is proved.
   
   Finally the last equation in the theorem follows by noting that 
   $$\Lambda_i\equiv \Lambda_i^{(2)}=1+\pi Z_i^{1}\ \ \textup{mod}\ \ \pi^2$$
   and using equations (\ref{omega}), (\ref{betabeta1}), (\ref{gamma}), (\ref{monster1}).
  \qed

\bigskip

{\it Proof of Theorem \ref{LCC}}.
Take, in Theorem \ref{theone}, $S=\mathcal A$, let $A_i$ and $B$ be the matrices $A_i^{(s)}$ and $B^{(s)}$ in Equation \ref{1955}, $P=\mathcal P=(x-1)$, and let  $L_{ij}^k:=L_{ij}^{k(s)}$. Then  the symmetry condition (\ref{uf22}) holds and therefore condition (\ref{uf2}) holds.
On the other hand note that
$$B\equiv (q^{(s)})^{(p^s)}\ \ \textup{mod}\ \ \mathcal P$$
hence 
$$L_{kij}\equiv \sum_l L_{ki}^l (q_{jl}^{(s)})^{p^s}=  L_{kij}^{(s)} \ \ \textup{mod}\ \ \mathcal P$$
and
$$\frac{1}{\pi}(B^{(s)}(\Lambda_i^{(s)}-1))_{kj}\equiv \sum_l (q^{(s)}_{kl})^{p^s}\Gamma_{ij}^{l(s)}=
\Gamma^{(s)\textup{LC}}_{ijk}\ \ \ \textup{mod}\ \ \mathcal P.
$$
We conclude by using  Lemma \ref{eajuns} and the last assertion of  Theorem \ref{theone}.\qed

\subsection{Arithmetic Chern connection}\label{acc} We now turn our attention to defining Chern connections. Set $G=\textup{GL}_{N/R_{\pi}}$, with $N$ not necessarily equal to $n$,
and recall that for $s\geq 1$ we have the trivial $\pi$-connection $((\delta_{1,0}^{(s)})^G,\ldots,(\delta_{n,0}^{(s)})^G)$  of degree $s$ on $G$. We consider its associated higher $\pi$-Frobenius lifts by $((\phi_{1,0}^{(s)})^G,\ldots,(\phi_{n,0}^{(s)})^G)$. We fix a metric $q^{(s)}\in G(R_{\pi})^{\textup{sym}}$, which we  view as an analogue of a metric on a vector bundle of rank $N$ on a manifold of dimension $n$. Generalizing \cite[Def. 4.25]{Bu17}, which corresponds to  the case $n=s=1$ and $\pi=p$, we introduce the following definition.

\begin{definition} Let $(\phi_i^{(s)})^G$ be a higher $\pi$-Frobenius lift of degree $s$  on $\mathcal A$ extending $\phi_i^{(s)}$ on $R_{\pi}$. We say that $(\phi_i^{(s)})^G$ is  {\bf $\cB_{q^{(s)}}$-symmetric} provided the following  diagram commutes
\begin{equation}
 \label{got-chern}
 \begin{array}{rcl}
\mathcal A& \stackrel{\cB_{q^{(s)}}}{\longrightarrow} & \mathcal A \otimes \mathcal A \\
\cB_{q^{(s)}}  \downarrow &\ &\downarrow (\phi_i^{(s)})^G \otimes (\phi_{i,0}^{(s)})^G\\
 \mathcal A \otimes \mathcal A & \stackrel{(\phi_{i,0}^{(s)})^G \otimes (\phi_i^{(s)})^G}{\longrightarrow} & \mathcal A\end{array}
\end{equation}
where $\cB_{q^{(s)}}$ is defined by $\cB_{q^{(s)}}(x) = x_1^tq^{(s)}x_2$ and $x$, $x_1:=x\otimes 1$, $x_2:=1\otimes x$ are the corresponding $N\times N$ matrices of indeterminates. 
\end{definition}

\begin{remark}\label{Bqcond}
As in \cite{BD} or \cite{Bu17} and with notation as in (\ref{1955})
the commutativity of the diagram (\ref{got-chern}) is equivalent to the equality
$$
A_i^{(s)} \Lambda_i^{(s)} = \Lambda_i^{(s)t}A_i^{(s)}.$$
\end{remark}

We give the PDE version of Chern connections in \cite{Bu17}.

\begin{theorem}\label{ch} For $q^{(s)} \in \textup{GL}_N(R_\pi)^{\textup{sym}}$ and $i \in \{1,\ldots,n\}$ there is a 
unique higher $\pi$-Frobenius lift $\phi_i^{(s)\textup{Ch}}$ of degree $s$ on $\mathcal A$ that is $\cB_{q^{(s)}}$-symmetric  and metric with respect to $q^{(s)}$; it is given by the formula
$$\phi_i^{(s)\textup{Ch}}(x) := x^{(p^s)} \{ (x^{(p^s)t}\phi_i^{(s)}(q^{(s)})x^{(p^s)})^{-1}(x^t q^{(s)} x)^{(p^s)}\}^{1/2}.$$ 
\end{theorem}

Here the matrix whose square root we are considering 
can be written in the form $1+\pi X$ with $X$ a matrix with entries in $\mathcal A$. For such a matrix the square root  is taken to be defined by the usual Newton binomial series; cf.  \cite[Rmk. 4.12]{Bu17}.

\begin{definition}\label{archcox}
The $\pi$-connection $\Delta^{(s)\textup{Ch}}=(\delta_1^{(s)\textup{Ch}},\ldots,\delta_n^{(s)\textup{Ch}})$  of degree $s$ 
defined by the higher Frobenius lifts $\phi_i^{(s)\textup{Ch}}$ is called the {\bf arithmetic Chern connection} of degree $s$ attached to $q^{(s)}$. \end{definition}

{\it Proof of Theorem \ref{ch}}.
This proof is entirely similar to the proof of \cite[Thm. 4.23]{Bu17}, by replacing $p$ with $\pi$
and working in degree $s$ as opposed to $s = 1$. We briefly recall the argument.  With notation as in (\ref{1955}) set $C_i^{(s)} = \frac{1}{\pi}(B^{(s)}-A_i^{(s)})$ and define 
$$\Lambda_i^{(s)} = \{(A_i^{(s)})^{-1}B^{(s)}\}^{1/2} = \{1 + \pi (A_i^{(s)})^{-1}C^{(s)}_i\}^{1/2}.$$
One has $\Lambda^{(s)}_i \equiv 1 \bmod \pi$. It is entirely formal to see that the $\mathcal B_{q^{(s)}}$-symmetry condition $A_i^{(s)} \Lambda_i^{(s)}= \Lambda_i^{(s)t} A_i^{(s)}$ holds; cf. Remark \ref{Bqcond}. One the other hand, the same algebraic manipulations in the proof of \cite[Thm. 4.23]{Bu17} demonstrate that the metric condition $\Lambda_i^{(s)t} A_i^{(s)} \Lambda_i^{(s)} = B^{(s)}$ holds; cf. Remark \ref{comua}. Existence now follows by noting that $\phi_i^{(s)\textup{Ch}}(x) = x^{(p^s)}\Lambda_i^{(s)}(x)$.  Uniqueness comes from the formal adaptation of \cite[Cor. 3.151 and 3.152]{Bu17}.
\qed

\begin{remark}
For $\pi = p$, $n = 1$, $s = 1$, and $N$ is arbitrary this is precisely the real Chern connection given in \cite[Thm. 4.23]{Bu17}. In particular, for $N = s=1$, this $\pi$-connection is given by $\phi^{\textup{Ch}}(x) = \left( \frac{q}{p} \right) q^{(p-1)/2} x^p$ where $\left( \frac{q}{p} \right)$ is the Legendre symbol. For more on the relation with the Legendre symbol we refer to Part 2 of this paper.
\end{remark}

 For the next statement consider $N = n$, $q^{(s)}\in \textup{GL}_n(R_{\pi})^{\textup{sym}}$,  and assume $L^{(s)}$ is a torsion symbol of the second kind so we may consider both the arithmetic Chern connection $\Delta^{\textup{(s)Ch}}:=(\delta_1^{\textup{(s)Ch}},\ldots,\delta_n^{\textup{(s)Ch}})$ and the Levi-Civita connection $\Delta^{\textup{(s)LC}}:=(\delta_1^{\textup{(s)\textup{LC}}},\ldots,\delta_n^{\textup{(s)\textup{LC}}})$ attached to our data. Denote by $(\Gamma_1^{\textup{(s)Ch}},\ldots,\Gamma_n^{\textup{(s)Ch}})$ the Christoffel symbols of the second kind for the Chern connection $\Delta^{\textup{(s)Ch}}$ and 
denote by $(\Gamma_1^{\textup{(s)LC}},\ldots,\Gamma_n^{\textup{(s)LC}})$ the Christoffel symbols of the second kind for the Levi-Civita connection $\Delta^{\textup{(s)LC}}$.
Also consider the associated Christoffel and torsion symbols of the first kind.

\begin{corollary}\label{tzinminte}
 We have the following congruences mod $\mathcal M$ in the ring $\mathcal A$:

\medskip

1)  $\Gamma_i^{\textup{(s)\textup{Ch}}} \equiv -\frac{1}{2} \delta_i^{(s)} q^{(s)} \cdot ((q^{(s)})^{(p^s)})^{-1}$.

\medskip

2) $\Gamma_{ijk}^{\textup{(s)Ch}} \equiv -\frac{1}{2} \delta_i^{(s)} q^{(s)}_{jk}.$

\medskip

3) $\Gamma_{ijk}^{\textup{(s)LC}} \equiv \frac{1}{2}\left(\Gamma_{ijk}^{\textup{(s)Ch}} + \Gamma_{jik}^{\textup{(s)Ch}} - \Gamma_{kij}^{\textup{(s)Ch}}\right)
+\frac{1}{2}
(L^{(s)}_{kij}(1)+L^{(s)}_{ijk}(1)-L^{(s)}_{jki}(1)).$

\end{corollary}

\

{\it Proof.}
The first congruence is entirely similar to \cite[Cor. 4.42, part (4.85)]{Bu17} and the second follows easily from the first.  The third congruence follows.\qed

\section{$\pi$-jet spaces and $\delta$-Lie algebra of $GL_n$}\label{deltalie}
In this section we give an intrinsic interpretation of our arithmetic PDE Riemannian  concepts along the lines of \cite{BD} and \cite[Sect. 4.1]{Bu17}. We are still adopting Conventions \ref{konvention1} and  \ref{konvention2}; hence  we fix $\pi\in \Pi$ and drop $\pi$ as an index, except in $R_{\pi}$. We also fix, in what follows, a family
$\Delta^{(s)}=(\delta_1^{(s)},\ldots,\delta_n^{(s)})$
of higher $\pi$-derivations of degree $s$ on $R_{\pi}$.

\subsection{Review of partial $\pi$-jet spaces}
In this subsection we  extend to arbitrary degree  the definition of  partial $\pi$-jet spaces \cite{BMP}; they are  universal initial objects in categories to be described presently which are similar to the $\pi = p$, $s=1$ case with one Frobenius lift \cite{Bu95}. See  \cite[Sec. 2.5]{BMP} for a more thorough treatment of partial $\pi$-jet spaces. 

\begin{definition}\label{proll} Let $s\in \mathbb N$.
Define the category ${\bf Prol}^{*(s)}$  as follows. An object of this category is  a countable family of $p$-adically complete $R_{\pi}$-algebras $S^* = (S^r)_{r \geq 0}$ equipped with the following data,
\begin{enumerate}
\item $R_{\pi}$-algebra homomorphisms $\varphi \colon S^r \to S^{r+1}$ for $r\geq 0$,
\item higher $\pi$-derivations $\delta_j^{(s)} \colon S^r \rightarrow S^{r+1}$ of degree $s$ for $1 \leq j \leq n$ and $r\geq 0$. 
\end{enumerate}

 \noindent  We require that $\delta_j^{(s)}$ be compatible with the $\pi$-derivations on $R_{\pi}$ and with $\varphi$, i.e.,  $\delta_j^{(s)} \circ \varphi = \varphi \circ \delta_j^{(s)}$. 
Morphisms are defined in a natural way. Denote the corresponding $\pi$-Frobenius lifts of degree $s$ by $\phi_j^{(s)} \colon S^r\rightarrow S^{r+1}$ which are therefore given by the rule $$\phi^{(s)}_j(x)=\varphi(x)^{p^s}+\pi\delta_j^{(s)}x.$$ 
\end{definition}

The objects of ${\bf Prol}^{*(s)}$ are called {\bf prolongation sequences} of degree $s$ over $R_{\pi}$. We sometimes identify elements  $a\in S^r$ with the elements  $\varphi(a)\in S^{r+1}$ if no confusion arises and sometimes write $S^*=(S^r,\varphi,\delta^{(s)}_1,\ldots,\delta^{(s)}_n)$.

\begin{remark} Recall the following obvious properties of this category. \
\begin{enumerate}
\item If $S$ is a $p$-adically complete partial $\delta$-ring over $R_{\pi}$ of degree $s$   whose $\pi$-derivations are compatible with those on $R_{\pi}$ then the sequence $S^*=(S^r)$ with $S^r=S$ has a natural structure of object of ${\bf Prol}^{*(s)}$ with $\varphi$ the identity and obvious $\delta^{(s)}_j$'s. The initial object in ${\bf Prol}^{*(s)}$ is the sequence $R^*_{\pi}=(R_{\pi}^r)$ with $R^r_{\pi}:=R_{\pi}$. 
\item If $S^*=(S^r,\varphi,\delta^{(s)}_1,\ldots,\delta^{(s)}_n)$ is an  object of ${\bf Prol}^{*(s)}$ then the inductive limit
$\varinjlim_{\varphi} S^r$ has a natural structure of  partial $\delta$-ring of degree $s$.
\end{enumerate}
\end{remark}

 We introduce some distinguished objects in ${\bf Prol}^{*(s)}$ which, in the $s=1$ case, play a central role in the theory \cite{BMP}. For $y = (y_1,\ldots,y_N)$ a tuple of indeterminates, consider the ring $R_{\pi}[y]$. Recalling the Notation  \ref{MMn}
 consider new indeterminates  denoted by $(\delta^{(s)})_\mu y_j$ for $\mu\in {\mathbb M}_n$,  $j\in \{1,\ldots,N\}$. For $r\geq 0$ the ring 
\begin{equation}
\label{mottby}
J^r(R_{\pi}[y]):=R_{\pi}[(\delta^{(s)})_\mu y_j\ |\ \mu\in {\mathbb M}_n^r, j\in \{1,\ldots,N\}]^{\widehat{\ }}\end{equation} is called the {\bf partial $\pi$-jet algebra} of $R_{\pi}[y]$ of order $r$ and degree $s$.  The sequence 
$$J^*(R_{\pi}[y]):=(J^r(R_{\pi}[y]))$$ has a unique structure of object in ${\bf Prol}^{*(s)}$ such that $$\delta^{(s)}_i (\delta^{(s)})_\mu y:=(\delta^{(s)})_ {i\mu}y$$ for all $i=1,\ldots,n$.

For $F\in J^*(R_{\pi}[y])$ one defines 
the evaluation map $F_{R_{\pi}} \colon R_{\pi}^N\rightarrow R_{\pi}$ by sending $(a_1,\ldots,a_N)\in
R_{\pi}^N$ to  the element $F_{R_{\pi}}(a_1,\ldots,a_N)\in R_{\pi}$ obtained from $F$ by replacing the indeterminates $(\delta^{(s)})_{\mu}y_j$  with the elements $(\delta^{(s)})_{\mu}a_j$. Note that the map
\begin{equation}
\label{e0}
J^r(R_{\pi}[y])\rightarrow \text{Fun}(R_{\pi}^N, R_{\pi}),\ F\mapsto F_{R_{\pi}}
\end{equation} is not injective in general. 

For every 
  $R_{\pi}$-algebra of finite type $A := R_{\pi}[y]/I$, $I$ an ideal, we define the rings
$$J^r(A):=J^r(R_{\pi}[y])/((\delta^{(s)})_{\mu}I\ |\ \mu\in {\mathbb M}_n^r).$$ This algebra is called the {\bf partial $\pi$-jet algebra} of $A$  of order $r$ and degree $s$. Even though it depends on $s$, we suppress the notation as it causes no confusion in what follows. Note that the sequence $J^*(A):=(J^r(A))$ has a natural structure of prolongation sequence i.e., it is an object of ${\bf Prol}^{*(s)}$ and enjoys the expected universal property that for every object $T^*$ of ${\bf Prol}^{*(s)}$ and every $R_{\pi}$-algebra map $u:A\rightarrow T^0$ there is a unique morphism  $J^*(A)\rightarrow T^*$ over $S^*$ in  ${\bf Prol}^{*(s)}$ compatible with $u$. See \cite[Prop. 2.22]{BMP} for the $s=1$ case, the  general case follows similarly.

The jet construction can be globalized in the natural way as in   \cite[Prop. 2.21]{BMP}. For every 
 scheme $X$ of finite type over $R_{\pi}$, the $p$-adic formal scheme 
$$J^r(X)=\bigcup \textup{Spf}(J^r(\mathcal O(U_i))),$$
is called the {\bf partial $\pi$-jet space} of $X$ of order $r$ and degree $s$, where 
$X=\bigcup U_i$ is an affine open cover and the implied gluing is well defined as explained in loc.cit. The elements of the ring $\mathcal O(J^r(X))$, identified with morphisms of $p$-adic formal schemes $J^r(X)\rightarrow \widehat{\mathbb A^1}$,
are called  {\bf arithmetic PDEs} on $X$ over $R_{\pi}$ of order (at most) $r$ and degree $s$.

\begin{remark}
In our discussion $\pi$ and $s$ were fixed so they were not included in the notation $J^r(A)$, $J^r(X)$, etc. If one allows $\pi$ and $s$  to vary then the various partial $\pi$-jet algebras/spaces are linked by various natural maps. The case when  $\pi$ is variable and $s$ is fixed   plays a role in  Section \ref{sec:OC} below. For the case when $\pi$ is fixed and $s$ is variable an example of such natural maps, implicitly related to our discussion in  Remark \ref{coburn}, is as follows. 
 Let $\phi^{(1)}_1,\ldots,\phi^{(1)}_n$ be higher $\pi$-Frobenius lifts of degree $1$ on $R_{\pi}$,  view $R_{\pi}$ as a partial $\delta$-ring of degree $1$ with respect to the higher $\pi$-derivations of degree $1$ attached to $\phi^{(1)}_1,\ldots,\phi^{(1)}_n$, and denote by $J^{m(1)}(A)$ the  partial $\pi$-jet algebra of order $m$ and degree $1$ of a `variable' finitely generated $R_{\pi}$-algebra $A$. On the other hand let $\nu_1,\ldots,\nu_n\in \mathbb M_n^{(s)}$ be some fixed words of some fixed length $s\in \mathbb N$, view $R_{\pi}$ as a partial $\delta$-ring of degree $s$ with respect to the higher $\pi$-derivations of degree $s$ attached to $\phi^{(s)}_1,\ldots,\phi^{(s)}_n$, where
 $\phi^{(s)}_j:=\phi^{(1)}_{\nu_j}$, and  denote by $J^{r(s)}(A)$ the partial $\pi$-jet algebra of order $r$ and degree $s$ of a `variable' finitely generated $R_{\pi}$-algebra $A$. Then for all $r\in \mathbb N$ there are  unique $R_{\pi}$-algebra 
homomorphism 
\begin{equation}
\label{astridd}
J^{r(s)}(A)\rightarrow J^{rs(1)}(A),\end{equation}
functorial in $A$,
with the property that for all $a\in A$, all $l\geq 1$, and all $\mu=i_1\ldots i_l\in \mathbb M_n^{(l)}$ the homomorphism (\ref{astridd}) sends
$$\phi^{(ls)}_{\mu}a\mapsto \phi^{(1)}_{\nu_{i_1}}\ldots \phi^{(1)}_{\nu_{i_l}} a.
$$
\end{remark}

\subsection{Group structure on partial $\pi$-jet spaces}\label{gs}
Recall that we considered the group scheme 
$G=\textup{GL}_{N/R_{\pi}}=\textup{Spec}(R_{\pi}[x,\det(x)^{-1}])$
 where $x$ is an $N\times N$ matrix of indeterminates.  One can consider then the  partial $\pi$-jet space $J^1(G)$ of order $1$ and degree $s$;   one can easily check that
$$J^1(G)\simeq \textup{Spf}(R_{\pi}[x,\delta^{(s)}_1 x,\ldots,\delta^{(s)}_n x,\det(x)^{-1}])^{\widehat{\ }}.$$
Consider the matrix of indeterminates $\delta^{(s)}x:=(\delta^{(s)}x_{jk})$  and define the $p$-adic formal scheme
$$\mathfrak g:=\textup{Spf}(R_{\pi}[\delta^{(s)} x])^{\widehat{\ }}.$$
Hence
$$\mathfrak g^n=\textup{Spf}(R_{\pi}[\delta_1^{(s)} x,\ldots,\delta^{(s)}_n x])^{\widehat{\ }}.$$
We will presently define a group structure on $\mathfrak g$ in the category of $p$-adic formal schemes but for now we ignore this. For every $p$-adically complete $R_{\pi}$-algebra $S$ the set of $S$-points (i.e., $\textup{Spf}(S)$-points)
$\mathfrak g(S)$ identifies with the set
$\textup{Mat}_n(S)$.
We have a canonical identification in the category of $p$-adic formal schemes 
$$J^1(G)\simeq \widehat{G}\times \mathfrak g^n.$$
For $S$ as above the $S$-points  of the formal scheme 
$\widehat{G}\times \mathfrak g$ are identified with tuples
$$(a_0,a_1,\ldots,a_n),\ \ a_0\in \textup{GL}_n(S),\ \ a_1,\ldots,a_n\in \textup{Mat}_n(S).$$
On the other hand by the universal property of partial $\pi$-jet spaces we have an induced structure of group (in the category of $p$-adic formal schemes) on $J^1(G)$. This induces a 
structure of group (in the category of $p$-adic formal schemes) on $\widehat{G}\times \mathfrak g$. As in \cite[p. 171]{Bu17} this group structure is defined on $S$-points by the multiplication rule:
$$(a_0,\ldots,a_n)(b_0,\ldots,b_n)=(c_0,\ldots,c_n)$$
where $c_0=a_0b_0$ and for every $i\in \{1,\ldots,n\}$ we have:
$$c_i=a_0^{(p^s)}b_i+ a_ib_0^{(p^s)}+\pi a_i b_i+\pi^{-1}(a_0^{(p^s)}b_0^{(p^s)}-(a_0b_0)^{(p^s)}).
$$
We have a natural identification of $\mathfrak g^n(S)$ with 
the set of $S$-points of the form  $$(1,a_1,\ldots,a_n)$$ of $\widehat{G}\times \mathfrak g^n$ and hence an identification of $\mathfrak g^n$
the  kernel of the group homomorphism
$J^1(G)\rightarrow \widehat{G}$. The induced group structure on $\mathfrak g^n$ 
is simply the $n$-fold product of the group structure on $\mathfrak g$ which, on
 $S$-points, is given by the operation
\begin{equation}
\label{zu1}
\mathfrak g(S)\times \mathfrak g(S)\rightarrow \mathfrak g(S),\ \ (a,b)\mapsto a +_{{\pi}}b:=a+b+\pi ab.\end{equation}
Note that the inverse in $\mathfrak g(S)$ is given by the map
\begin{equation}
\label{zu2}
\mathfrak g(S)\rightarrow \mathfrak g(S),\ \ a\mapsto -a+\pi a^2-\pi^2 a^3+\ldots.\end{equation}
If we set
\begin{equation}
\label{G1}
G^1(R_{\pi}):=\textup{ker}(G(R_{\pi})\rightarrow G(k))\end{equation}
then we have a group isomorphism
\begin{equation}
\label{gG}
\mathfrak g(R_{\pi})\rightarrow G^1(R_{\pi}),\ \ a\mapsto 1+\pi\cdot a.\end{equation}
The group $\mathfrak g^n$ can be viewed as a `$p$-adic formal analogue' of the formal group of $\textup{GL}_N$.
By analogy with the classical theory we refer to the group $\mathfrak g^n$ as the 
{\bf $\delta$-Lie algebra} of $G$ of degree $s$, cf., \cite[Def. 3.132]{Bu17}. Note, however,  that, as a group,
$\mathfrak g(S)$, and hence $\mathfrak g^n(S)$, is non-abelian for $N\geq 2$. Also note that
the various groups $\mathfrak g^n$ as $s$ varies in $\mathbb N$ are canonically isomorphic.

\subsection{Logarithmic derivative and Christoffel symbols}\label{ld}

 By the universal property of partial $\pi$-jet spaces of degree $s$ for every $p$-adically complete partial $\delta$-ring $S$ of degree $s$ we have an induced map of sets
$$D:G(S)\rightarrow J^1(G)(S)\simeq G(S)\times \mathfrak g^n(S),\ \ D (a)=(a,\delta^{(s)}_1 a,\ldots, \delta^{(s)}_n a).$$
Assume now 
that we are given a  $\pi$-connection $((\delta_1^{(s)})^G,\ldots,(\delta_n^{(s)})^G)$ of degree $s$ on $G$. Taking $S=\mathcal A$  and $a=x^t$ we have 
$$D (x^t)=(x^t,(\delta_1^{(s)})^Gx^t,\ldots,(\delta_n^{(s)} )^G x^t).$$
Let $(\Gamma^{(s)}_1,\ldots,\Gamma^{(s)}_n)$ be the Christoffel symbols of the second kind attached to our  $\pi$-connection.
If we consider in addition the trivial $\pi$-connection $((\delta_{1,0}^{(s)})^G,\ldots,(\delta_{n,0}^{(s)})^G)$ of degree $s$ on $G$ (with $(\delta_{i,0}^{(s)})^G x=0$) and we consider the corresponding map $D_0$ then, using the group law on the partial $\pi$-jet space, we have the following computation:
$$\begin{array}{rcl}
(D(x^t))(D_0(x^t))^{-1} & = & (1,(\delta_1^{(s)})^G x^t \cdot (x^{(p^s)t})^{-1},\ldots,(\delta_n^{(s)})^G x^t \cdot (x^{(p^s)t})^{-1})\\
\ & \ & \ \\
\ & = & (1,\Gamma^{(s)}_1,\ldots,\Gamma^{(s)}_n).\end{array}$$
So we can interpret the Christoffel symbols of the second kind as a tuple
$$(\Gamma^{(s)}_1,\ldots,\Gamma^{(s)}_n) \in \mathfrak g(S)^n$$
obtained by taking the `difference' between $D$ and $D_0$ evaluated at $x^t$. This difference is an analogue of the classical `logarithmic derivative' map or 
`Maurer-Cartan equations' of a Lie group as explained in \cite[loc.cit]{Bu17}.

 \section{$\delta$-overconvergence}\label{sec:OC}

 We now consider overconvergence aspects of the theory.  This involves 
 keeping $s$ fixed and allowing  $\pi\in \Pi$ to vary. 
 We therefore return  to using more subscripts $p$, $\pi$, etc., to make clear over which rings various objects are defined; in other words we do not adopt here Convention \ref{konvention1}.
 Instead,  throughout our discussion here, we fix a tuple of higher Frobenius automorphisms $\Phi^{(s)}=(\phi^{(s)}_1,\ldots,\phi^{(s)}_n)$ on $K^{\textup{alg}}$ of degree $s$ and for every $\pi\in \Pi$ we consider the induced  tuple of higher $\pi$-Frobenius lifts $\Phi^{(s)}_{\pi}=(\phi^{(s)}_{\pi,1},\ldots,\phi^{(s)}_{\pi,n})$ of degree $s$ on  $R_{\pi}$ and the 
 corresponding tuple of higher $\pi$-derivations  $\Delta^{(s)}_{\pi}=(\delta^{(s)}_{\pi,1},\ldots,\delta^{(s)}_{\pi,n})$ of degree $s$ on $R_{\pi}$. 
 We denote by ${\bf Prol}^{*(s)}_{\pi}$  the  category of prolongation sequences of degree $s$ corresponding to $\pi$ and $\Delta^{(s)}_{\pi}$ and  denote by $J^r_{\pi}=J^{r(s)}_{\pi}$ the  partial $\pi$-jet space/algebra functors of order $r$ and degree $s$ corresponding to $\pi$ and $\Delta^{(s)}_{\pi}$; cf. Section \ref{deltalie}. We denote by ${\bf Prol}^{(s)}_{\pi}$ the subcategory of all $S^*$ for which $S^r$ are Noetherian and flat over $R_{\pi}$. If $A$ is smooth over $R_{\pi}$ then $J^*_{\pi}(A):=(J^r_{\pi}(A))$ is an object of ${\bf Prol}^{(s)}_{\pi}$; cf. \cite[Prop. 2.22]{BMP}.

 We first review $\delta$-overconvergence extending the discussion in  \cite[Sec. 2.6]{BMP} where the case $s=1$ was considered. This builds on ideas initiated in \cite{BS11,BM20}.

Exactly as in \cite[Sect. 4.1]{BM20}
for every $\pi' | \pi$ and  every object $S^*$ in ${\bf Prol}^{(s)}_{\pi}$ the sequence  
$$S^* \otimes_{R_{\pi}} R_{\pi'} := (S^r \otimes_{R_{\pi}} R_{\pi'})_{r \geq 0}$$ is naturally an object of ${\bf Prol}^{(s)}_{\pi'}$.  For a scheme $X$ of finite type over $R_{\pi}$ we set $X_{\pi'}:=X\otimes_{R_{\pi}} R_{\pi'}$.
Clearly $J^0_{\pi'}(X_{\pi'})=\widehat{X_{\pi'}}$. Note also that
$J^r_{\pi'}(X_{\pi'})$ only depends on $r,\pi',X$, and 
 $\Phi^{(s)}_{\pi}$. 

When $A$ is a smooth $R_{\pi}$-algebra, for all $\pi''|\pi'|\pi$ there  are natural homomorphisms
\begin{equation}
\label{e1}
\iota_{\pi'',\pi'}:
J^r_{\pi''}(A)\rightarrow J^r_{\pi'}(A) \otimes_{R_{\pi'}} R_{\pi''}\end{equation}
such that  the homomorphism
\begin{equation}
\label{e2}
\iota_{\pi'',\pi}:J^r_{\pi''}(A)\rightarrow 
J^r(A) \otimes_{R_{\pi}} R_{\pi''}
\end{equation}
equals the composition
\begin{equation}
\label{e3}
J^r_{\pi''}(A)\stackrel{\iota_{\pi'',\pi'}}{\longrightarrow} J^r_{\pi'}(A) \otimes_{R_{\pi'}} R_{\pi''}\stackrel{\iota_{\pi',\pi}\otimes 1}{\longrightarrow} 
(J^r(A) \otimes_{R_{\pi}} R_{\pi'})\otimes_{R_{\pi'}}R_{\pi''},
\end{equation}
where the targets of the maps (\ref{e2}) and (\ref{e3}) are naturally identified. Moreover the homomorphisms \eqref{e1} are injective; this can be checked exactly as in  \cite[Prop. 2.25]{BMP}.

The notion of $\delta$-overconvergence was first introduced in \cite{BS11} and then exploited in \cite{BM20}. The following version was introduced in  \cite[Def. 2.27]{BMP} for $s=1$.

\begin{definition} \label{overconvergence}
For $A$  a smooth $R_{\pi}$-algebra, an element  $f_{\pi}\in J^r_{\pi}(A)$ is called {\bf totally $\delta$-overconvergent} provided for all $\pi'|\pi$ there is an integer $N\geq 0$ so that $p^N f_{\pi}\otimes 1$ is in the image of 
the map 
\begin{equation}
\label{e2t}\iota_{\pi',\pi}:J^r_{\pi'}(A)\rightarrow 
J^r_{\pi}(A) \otimes_{R_{\pi}} R_{\pi'}.\end{equation} 
 \end{definition}
 
 Like in  \cite{BMP} we denote by $J^r_{\pi}(A)^!$ the $R_{\pi}$-algebra of  totally $\delta$-overconvergent elements in $J^r_{\pi}(A)$; then,  for $X:=\textup{Spec}(A)$ and every $f\in J^r_{\pi}(A)^!$ we have a naturally induced map of sets
 $$f^{\textup{alg}}:X(R^{\textup{alg}})\rightarrow K^{\textup{alg}}.$$

\bigskip

We now introduce a variation of the above concepts.

\begin{definition}
Let $S$ be a smooth $R$-algebra.
A matrix 
$$\Theta_p=(\Theta_{p,ij})\in \textup{Mat}_n(J^r_{p}(S))[x,\det(x)^{-1}]^{\widehat{\ }})$$
 will be called   {\bf  strictly $\delta$-overconvergent} if for all $\pi\in \Pi$ the images 
 of the entries $\Theta_{p,ij}$ via the (injective) homomorphism
 \begin{equation}
 J^r_{p}(S)[x,\det(x)^{-1}]^{\widehat{\ }}\rightarrow (J^r_{p}(S)\otimes_R R_{\pi})[x,\det(x)^{-1}]^{\widehat{\ }}\end{equation}
  belong to the image of the (injective) ring homomorphism
 \begin{equation}
 \label{irh}
 J^r_{\pi}(S)[x,\det(x)^{-1}]^{\widehat{\ }}\rightarrow (J^r_{p}(S)\otimes_R R_{\pi})[x,\det(x)^{-1}]^{\widehat{\ }}.\end{equation}\end{definition}
 
 We shall apply this to the case
 $$S=R[Q,\det(Q)^{-1}]$$
  where $Q:=(Q_{ij})$ is a symmetric matrix with indeterminate coefficients.

 \begin{remark}\ 
 
 1) The injectivity of the maps above follows from the analogue (for $s$ arbitrary) of  \cite[Prop. 2.25]{BMP} plus
 the fact that for every Noetherian $p$-adically complete flat $R_{\pi}$-algebra  $A$ the natural map 
 $$A[x,\det(x)^{-1}]^{\widehat{\ }}\rightarrow A[[x-1]]$$
  is injective, as one can see by reducing this map modulo $\pi$. 
 
 2) If $\Theta_p$ above is strictly $\delta$-overconvergent
 then $\Theta_{p,ij}$, viewed as elements of $J^r_{p}(S[x,\det(x)^{-1}])$, are totally $\delta$-overconvergent in the sense of Definition \ref{overconvergence}.

 3) 
 For every strictly $\delta$-overconvergent  $\Theta_p$ as above we have well defined maps
 $$\Theta_p^{\textup{alg}}:\textup{Mat}_n(R_{\pi})^{\textup{sym}}\rightarrow \textup{Mat}_n(R^{\textup{alg}}[x,\det(x)^{-1}]^{\widehat{\ }}),\ \ q\mapsto \Theta_p^{\textup{alg}}(q).$$
 If $q$ has coefficients in $R_{\pi}$ and $\Theta_{p,ij}$ are images of $\Theta_{\pi,ij}\in
 J^r_{\pi}(S)[x,\det(x)^{-1}]^{\widehat{\ }}$ then one defines 
 $\Theta_p^{\textup{alg}}(q)$ to be the matrix with entries the series in $R^{\textup{alg}}[x,\det(x)^{-1}]^{\widehat{\ }}$ obtained by replacing $\delta^{(s)}_{\pi,\mu} Q$ in $\Theta_{\pi,ij}$ by $\delta^{(s)}_{\pi,\mu}q$. This  construction is compatible (in an obvious sense) with the natural maps
 $$\Theta_{p,ij}^{\textup{alg}}:\textup{Mat}_n(R_{\pi})^{\textup{sym}}\times \textup{GL}_n(R_{\pi})\rightarrow K^{\textup{alg}}$$ 
 obtained by viewing $\Theta_{p,ij}$ as an element of $J^r_{p}(S[x,\det(x)^{-1}])^!$.
 \end{remark}
 
 \begin{theorem}\label{strictt}
 Fix a torsion symbol $L^{(s)}_p=(L_p^{1(s)},\ldots,L_p^{n(s)})$ of the second kind with coefficients in $R$. Then there exist strictly $\delta$-overconvergent matrices
$$\Lambda^{(s){\bf LC}}_{p,1},\ldots, \Lambda^{(s){\bf LC}}_{p,n}\in \textup{Mat}_n(J^1_{p}(S))[x,\det(x)^{-1}]^{\widehat{\ }})$$
satisfying the following property. For all $\pi\in \Pi$ and all $q^{(s)}\in \textup{Mat}_n(R_{\pi})^{\textup{sym}}$
if
 $\Phi^{(s)\textup{LC}}=(\phi_{\pi,1}^{(s)\textup{LC}},\ldots,\phi_{\pi,n}^{(s)\textup{LC}})$ is the 
 tuple of higher $\pi$-Frobenius lifts of degree $s$ corresponding  to the arithmetic Levi-Civita connection of degree $s$ attached to 
$q^{(s)}$  and  to the torsion symbol $L_{\pi}^{(s)}:=\frac{p}{\pi}\cdot L_p^{(s)}$
    then for all $i\in \{1,\ldots,n\}$ we have
$$\phi^{(s)\textup{LC}}_{\pi,i} (x)=x^{(p^s)}\cdot (\Lambda_{p,i}^{(s){\bf LC}})^{\textup{alg}}(q^{(s)}).$$
\end{theorem}

{\it Proof}. 
For $\pi\in \Pi$ consider the matrices
 $$A_{\pi,i}^{(s){\bf LC}}:=x^{(p^s)t}\phi^{(s)}_{\pi,i}(Q) x^{(p^s)},\ \ B^{(s){\bf LC}}=(x^tQx)^{(p^s)}.$$
 The matrix $ B^{(s){\bf LC}}$ does not depend on $\pi\in \Pi$; similarly the images of the matrices $A_{\pi,i}^{(s){\bf LC}}$ via the homomorphism (\ref{irh}) are all equal and we denote by $A_i^{(s){\bf LC}}$ this common image.
By Theorem \ref{theone} applied to $L_{\pi}^{k(s)}=(L_{\pi,ij}^{k(s)})$ there exist 
a unique tuple $\Lambda_{\pi}^{(s){\bf LC}}$ of matrices
$$\Lambda^{(s){\bf LC}}_{\pi,1},\ldots, \Lambda^{(s){\bf LC}}_{\pi,n}\in \textup{Mat}_n(J^1_{\pi}(S)[x,\det(x)^{-1}]^{\widehat{\ }})$$
such that
\begin{equation}
\label{lpi}
\begin{array}{rcl}
(\Lambda^{(s){\bf LC}}_{\pi,i})^tA_i^{(s){\bf LC}} \Lambda_{\pi,i}^{(s){\bf LC}} & = & B^{(s){\bf LC}},\\
\ & \ & \ \\
(\Lambda_{\pi,i}^{(s){\bf LC}}-1)_{kj}  - (\Lambda_{\pi,j}^{(s){\bf LC}}-1)_{ki}
& = & \pi \cdot L_{\pi,ij}^{k(s)}(\Lambda_{\pi}^{(s){\bf LC}})\\
\ & \ & \ \\
\ & = & 
p \cdot L_{p,ij}^{k(s)}(\Lambda_{\pi}^{(s){\bf LC}}).\end{array}\end{equation}
By the uniqueness property in Theorem \ref{theone} we get that
$$\Lambda_{\pi,i}^{{\bf LC}}=\Lambda_{p,i}^{{\bf LC}}$$
for all $i\in \{1,\ldots,n\}$ and all $\pi\in \Pi$. In particular we get that 
$\Lambda_{p,i}^{(s){\bf LC}}$ are strictly $\delta$-overconvergent. 
Specializing $\delta^{(s)}_{\pi}Q\mapsto \delta^{(s)}_{\pi}q^{(s)}$ in Equations \ref{lpi} we get
\begin{equation}
\label{lpii}
\begin{array}{rcl}
((\Lambda^{(s){\bf LC}}_{p,i})^{\textup{alg}}(q^{(s)}))^t((A_i^{(s){\bf LC}})^{\textup{alg}}(q^{(s)})) ((\Lambda_{p,i}^{(s){\bf LC}})^{\textup{alg}}(q^{(s)})) & = & (B^{(s){\bf LC}})^{\textup{alg}}(q^{(s)})\\
\ & \ & \ \\
((\Lambda_{p,i}^{(s){\bf LC}})^{\textup{alg}}(q^{(s)})-1)_{kj} - ((\Lambda_{p,j}^{(s){\bf LC}})^{\textup{alg}}(q^{(s)})-1)_{ki} & = & p L_{p,ij}^k(\Lambda_p^{(s){\bf LC}}(q^{(s)})).\end{array}\end{equation}
On the other hand 
if we write
$$A_i^{(s)}=(A_i^{(s){\bf LC}})^{\textup{alg}}(q^{(s)})=x^{(p^s)t}\phi^{(s)}_{\pi,i}(q^{(s)}) x^{(p^s)},$$
$$B^{(s)}=(B^{(s){\bf LC}})^{\textup{alg}}(q^{(s)})=(x^tq^{(s)}x)^{(p^s)},$$
$$\phi_{\pi,i}^{(s)\textup{LC}} x=x^{(p^s)}\cdot \Lambda^{(s)}_{\pi,i}$$ then by the definition of the arithmetic Levi-Civita connection attached  to $q^{(s)}$ and $L_{\pi}^{(s)}$
we  have
\begin{equation}
\label{lpiii}
\begin{array}{rclll}
\Lambda^{(s)t}_{\pi,i}A_i^{(s)} \Lambda_{\pi,i}^{(s)} & = & B^{(s)}, &\  & \ \\
\ & \ & \ & \ & \ \\
(\Lambda^{(s)}_{\pi,i}-1)_{kj} - (\Lambda^{(s)}_{\pi,j}-1)_{ki} & = & \pi\cdot 
L^{k(s)}_{\pi, ij}(\Lambda^{(s)}_{\pi}) & = &
p \cdot L^{k(s)}_{p,ij}(\Lambda^{(s)}_{\pi}).\end{array}\end{equation}
By the uniqueness in the definition of the arithmetic Levi-Civita connection the equations
\ref{lpii} and \ref{lpiii} imply that 
$$\Lambda^{(s)}_{\pi,i}=(\Lambda_{p,i}^{(s){\bf LC}})^{\textup{alg}}(q^{(s)}).$$
\qed

\bigskip

A similar picture holds for the Chern connection. Indeed Theorem \ref{ch} immediately implies the following:

 \begin{theorem}\label{javasp}
  Then there exist strictly $\delta$-overconvergent matrices
$$\Lambda^{(s){\bf Ch}}_{p,1},\ldots, \Lambda^{(s){\bf Ch}}_{p,n}\in \textup{Mat}_n(J^1_{p}(S))[x,\det(x)^{-1}]^{\widehat{\ }})$$
satisfying the following property. For all $\pi\in \Pi$ and all $q^{(s)}\in \textup{Mat}_n(R_{\pi})^{\textup{sym}}$
if
 $\Phi^{(s)\textup{Ch}}=(\phi_{\pi,1}^{(s)\textup{Ch}},\ldots,\phi_{\pi,n}^{(s)\textup{Ch}})$ is the 
 tuple of higher $\pi$-Frobenius lifts of degree $s$ corresponding  to the arithmetic Chern connection of degree $s$ attached to 
$q^{(s)}$ 
    then for all $i\in \{1,\ldots,n\}$ we have
$$\phi^{(s)\textup{Ch}}_{\pi,i} (x)=x^{(p^s)}\cdot (\Lambda_{p,i}^{(s){\bf Ch}})^{\textup{alg}}(q^{(s)}).$$
\end{theorem}


\begin{thebibliography}{10}

\bibitem[BPS]{BPS} A. Bertapelle, E. Previato, A. Saha, {\it Arithmetic jet spaces},  arXiv:2003.12269, 2020.



\bibitem[Bu95]{Bu95} A. Buium, {\it Differential characters of Abelian varieties over $p-$adic fields}, Invent. Math., 122,  (1995), 309--340.
     
\bibitem[Bu96]{Bu96} A. Buium, {\it  Geometry of p-jets},  Duke Math. J., 82, 2, (1996), 349--367.
   
\bibitem[Bu97]{Bu97} A. Buium, {\it Differential characters and characteristic polynomial of Frobenius}, Crelle J., 485 (1997), 209--219.


\bibitem[Bu05]{Bu05} A. Buium, {\it Arithmetic Differential Equations},  Math. Surveys and Monographs, 118, AMS (2005).

\bibitem[Bu17]{Bu17} A. Buium, {\it Foundations of Arithmetic Differential Geometry}, Math. Surveys and Monographs, 222, AMS, 2017.

\bibitem[Bu19]{Bu19} A. Buium, {\it Arithmetic Levi-Civita connection},   Sel. Math. New Ser. 25, 12 (2019), 79 pages.

\bibitem[Bu20]{Bu20} A. Buium, {\it Arithmetic analogues of Hamiltonian systems}, in: Integrable Systems and Algebraic Geometry, 2, London Math. Soc., Lecture Notes Series 459, Cambridge University Press, 2020, 13--40.

\bibitem[BD16]{BD} A. Buium, T. Dupuy,  {\it Arithmetic differential equations on $GL_n$, II: arithmetic Lie-Cartan theory}, Selecta Math.  22, 2, (2016), 447--528.



\bibitem[BP09]{BP09} A. Buium, B. Poonen, {\it Independence of points on elliptic curves arising from special points on modular and Shimura curves, II: local results}, Compositio Math., 145 (2009), 566--602.
 
 \bibitem[BM15]{BMan} A. Buium, Yu. I. Manin, {\it  Arithmetic differential equations of Painlev\'{e} VI type},  in: Arithmetic and Geometry, London Mathematical Society Lecture
Note Series: 420,  L. Dieulefait, G. Faltings, D. R. Heath-Brown, Yu. V. Manin, B. Z. Moroz and J.-P. Wintenberger (eds),  Cambridge University Press, (2015), 114--138.


\bibitem[BM20]{BM20} A. Buium, L.E. Miller, {\it Solutions to arithmetic differential equations in algebraically closed fields}, Advances in Math. 375 (2020), 107342, 47pp.

\bibitem[BMa]{BMP} A. Buium, L.E. Miller, {\it Purely arithmetic PDE'S over a  $p$-adic field, I: $\delta$-characters and $\delta$-modular forms},  arXiv:2103.16627.

\bibitem[BMb]{BMadg2} A. Buium, L. E. Miller, {\it Arithmetic differential geometry in the PDE setting, II: curvature}, in preparation.

\bibitem[BP17]{BP17} A. Buium,  E. Previato, {\it Arithmetic Euler top},   J. Number Theory, 173, (2017), 37--63.

\bibitem[BP18]{BP18} A. Buium, E. Previato, {\it The Euler top and canonical lifts}, J. Number Theory, 190, (2018), 156--168.

\bibitem[BS11]{BS11} A. Buium, A. Saha,  {\it Differential overconvergence}, in: Algebraic methods in dynamical systems, Banach Center Publications, 94, (2011), 99--129.

\bibitem[BS09a]{BS09a} A. Buium, S.R. Simanca,  {\it  Arithmetic differential equations in several variables},   Annales Inst. Fourier, 59, 7 (2009), 2685--2708 (volume dedicated to B. Malgrange).

\bibitem[BS09b]{BS09b} A. Buium, S.R.  Simanca,  {\it  Arithmetic Laplacians},  Advances in Math.  220 (2009), 246--277.

\bibitem[BS10a]{BS10a} A. Buium, S.R. Simanca,  {\it Arithmetic partial differential equations I}, Advances in Math., 225, (2010), 689--783.  	

\bibitem[BS10b]{BS10b} A. Buium, S.R. Simanca, {\it Arithmetic partial differential equations II}, Advances in Math., 225, (2010), 1308--1340.	

\bibitem[Joy85]{J85} A. Joyal,  \textit{$\d-$anneaux et vecteurs de Witt}, C.R. Acad. Sci. Canada, VII, 3, (1985), 177--182.

\bibitem[Ked10]{Ked} K. Kedlaya, {$p$-adic differential equations}, Cambridge studies in advanced mathematics 125, CUP, 2010.

\bibitem[Ko84]{Ko84} N. Koblitz, {\it $p$-adic Numbers, $p$-adic analysis, and Zeta Functions}, Grad. Texts in Math. 58, Springer, 1984.

\bibitem[Man63]{Man63} Yu. I. Manin, {\it Rational points on algebraic curves over function fields}, Izv. Acad. Nauk, USSR, 27 (1963), 1395--1440.

\bibitem[Man13]{Man13} Yu. I. Manin, {\it Numbers as functions}, P-adic Numbers, Ultrametric Analysis and Applications, 5, 4, 2013, 313--325

\bibitem[Neu80]{N80} J. Neukirch, {\it Class Field Theory}, Springer, 1980.

\bibitem[vNeu53]{VN} J. von Neumann, {\it A certain zero-sum two-person game equivalent to an optimal assignment problem}, Ann. Math. Studies 28, 5-12, 1953.



\bibitem[Ser79]{Se79} J.-P. Serre, {\it Local Fields}, Springer, 1979.

\bibitem[VVZ94]{VVZ} V. S. Vladimirov, I. V. Volovich and E. I. Zelenov, {\it p-Adic Analysis and Mathematical Physics},  World Scientific, Singapore, 1994.

\bibitem[Was82]{Wa} L. Washington, {\it Introduction to cyclotomic fields}, Grad. Texts in Math. 83, Springer, 1982.

\end{thebibliography}
\end{document}